\documentclass[preprint,12pt,3p]{elsarticle}

\usepackage{amssymb}
\usepackage{amsmath}
\usepackage{yhmath}
\usepackage{amsthm}

\usepackage{url}
\usepackage{comment}
\usepackage{xcolor}
\usepackage{bigints}

\newtheorem{thm}{Theorem}
\newdefinition{lem}{Lemma}
\newdefinition{rmk}{Remark}
\newdefinition{definition}{Definition}
\newdefinition{cor}{Corollary}
\newproof{pf}{Proof}
\newcommand{\dd}{ {\rm d} }

\newcommand{\pv}{\operatorname{p.\!v.}}
\newcommand{\fp}{\operatorname{p.\!f.}}

\newcommand{\bn}{{\mathbf{n}}}
\newcommand{\btau}{{\boldsymbol{\tau}}}
\newcommand{\bgamma}{{\boldsymbol{\gamma}}}

\newcommand{\bx}{{\mathbf{x}}}
\newcommand{\by}{{\mathbf{y}}}
\newcommand{\br}{{\mathbf{r}}}
\newcommand{\nx}{{\bn(\bx)}}
\newcommand{\ny}{{\bn(\by)}}
\newcommand{\taux}{{\btau(\bx)}}
\newcommand{\tauy}{{\btau(\by)}}

\newcommand{\cK}{{\mathcal{K}}}

\DeclareMathOperator{\sign}{sign}

\journal{Journal of Computational Physics}

\begin{document}

\begin{frontmatter}

\title{Boundary Integral Formulations for Flexural Wave Scattering in Thin Plates}

\affiliation[1]{organization={Committee on Computational and Applied Mathematics, University of Chicago},
addressline={5747 S. Ellis Avenue},
postcode={60637},
city={Chicago, IL},
country={USA}}

\affiliation[2]{organization={Department of Mathematical Sciences, New Jersey Institute of Technology},
addressline={Cullimore Hall},
postcode={07102},
city={Newark, NJ},
country={USA}}

\author[1]{Peter Nekrasov\corref{cor1}}
\ead{pn3@uchicago.edu}

\author[1]{{Zhaosen} Su}

\author[2]{Travis Askham}
 
\author[1]{Jeremy G. Hoskins}

\cortext[cor1]{Corresponding author.  }

\begin{abstract}
In this paper, we develop second kind integral formulations  for flexural wave scattering problems involving the clamped, supported, and free plate boundary conditions. While the clamped plate problem can be solved with layer potentials developed for the biharmonic equation, the free plate problem is more difficult due to the order and complexity of the boundary conditions. In this work, we describe a representation for the free plate problem that uses the Hilbert transform to cancel singularities of certain layer potentials, ultimately leading to a Fredholm integral equation of the second kind. Additionally, for the supported plate problem, we improve on an existing representation to obtain a second kind integral equation formulation.  With these representations it is possible to solve flexural wave scattering problems with high-order-accurate methods, examine the far field patterns of scattering objects, and solve large problems involving multiple scatterers.
\end{abstract}

\begin{highlights}
\item This paper presents second kind integral formulations for three common boundary conditions for the flexural wave equation -- the clamped, {supported, and free} plates.
\item  An original integral representation for the free plate is developed by using the Hilbert transform to eliminate singularities resulting from higher order boundary conditions. 
\item These methods enable fast and accurate solutions for flexural wave scattering problems, which have a wide range of applicability in geophysics and engineering.
\end{highlights}


\end{frontmatter}

\section{Background}

{
Scientific interest in the vibrations of free plates began with Galileo and Hooke's observations of nodal patterns and culminated in Ernst Chladni’s famous demonstration at the Paris Academy in 1808. For the relevant history of these developments, see \cite{lindsay1966story,mclaughlin1998good}. Today, the theory of plates is frequently applied to the modeling and observation of elastic waves, which have been implicated in the collapses of both bridges \cite{Drabek2003} and ice shelves \cite{massom18}. These waves, known as flexural waves, pose a significant challenge to the integrity of both natural and built environments, making them a critical area of study. 
}

Flexural waves are typically modeled by the thin plate approximation, where the vertical displacement of {the plate is assumed to be} a function of two horizontal dimensions. {We will describe the equations for exterior domains; the interior
case is analogous. Let $\Omega = \bigcup_{i=1}^M \Omega_i \subset \mathbb{R}^2$ be the union of a finite collection of bounded,
simply connected domains, 
$\Omega_i$, with smooth boundaries and disjoint closures, and let $E = \mathbb{R}^2\setminus \overline{\Omega}$ denote the exterior
of this collection. The $\Omega_i$ represent holes, inclusions, or other obstacles in a plate of infinite
extent. For $\bx :=(x,y) \in E$,} the out-of-plane displacement $u(\bx )$ is described by the following fourth-order partial differential equation (PDE) {in non-dimensional form}:
\begin{align}
    \Delta^2 u - k^4 u = 0 \, , \qquad {\text{ in } E} \, , \label{flexural}
\end{align}
where $\Delta^2 := \displaystyle \frac{\partial^4 }{\partial x^4} + 2 \frac{\partial^4}{\partial x^2 \partial y^2} + \frac{\partial^4}{\partial y^4}$ is the biharmonic operator and $k$ is the {non-dimensionalized} wavenumber, {which depends on the frequency, plate thickness, elastic moduli, and characteristic lengthscale}. The first term in the equation {relates to} a bending force, while the second term {relates to} the inertia {of the plate after} assuming that the solution is harmonic in time {\cite[\S 25]{landau59}}. 
{For interior problems,} \eqref{flexural} can also be thought of as an eigenvalue problem corresponding to the biharmonic operator {\cite{kato1957estimation}}.

Since equation \eqref{flexural} is fourth-order, two boundary conditions are required to fully determine {the boundary value} problem. {In the present work, we are concerned with three sets of boundary conditions that frequently arise in thin plate theory~\cite{landau59, kato1957estimation, timoshenko1959theory, ritz09, Friedrichs1928}: the clamped,
supported, and free plate conditions. The boundary conditions for the clamped plate problem, also known as the Dirichlet problem, take the form:
\begin{align}
    \begin{cases}
    \displaystyle u = f_1 \, , \qquad &\text{ on } \partial \Omega \, ,  \\
    \displaystyle \frac{\partial u}{\partial n}  = f_2 \, , \qquad &\text{ on } \partial \Omega \, , 
    \end{cases} \label{clampedbcs}
\end{align} 
where $\partial/\partial n$ denotes the normal derivative. The first condition represents the displacement of the plate while the second condition represents the slope {along the normal vector, which we always take to be outward facing}.
These boundary conditions are common in applications in which an edge of the plate is fixed at the boundary~\cite{lakitosh2012analysis,kapania2000static,sergienko13,Nekrasov_MacAyeal_2023}.

The boundary conditions for the supported plate problem, also known as the hinged plate problem,
take the form:
\begin{align}
\begin{cases}
    u = f_1 \, , &\qquad \text{ on } \partial \Omega \, ,  \\
    \displaystyle \nu \Delta u + (1-\nu) \frac{\partial^2 u}{\partial n^2}= f_2 \, , &\qquad \text{ on } \partial \Omega \, , 
\end{cases} \label{supportedbcs}
\end{align}
where $\nu$ is the Poisson's ratio of the plate, {which takes on values between $-1$ and $1/2$ \cite{landau59}}. Again, the first condition corresponds to a displacement while the second condition corresponds to a bending moment along the boundary.
The supported plate boundary conditions are common in applications where an edge of the
plate is supported from below~\cite{seide1958,zhao2002plate,MEYLAND2021}.

Finally, the boundary conditions for the free plate problem, or the Neumann problem,
take the form: 
\begin{align}
\begin{cases}
    \displaystyle \nu \Delta u + (1-\nu) \frac{\partial^2 u}{\partial n^2}= f_1 \, ,  &\qquad \text{ on } \partial \Omega \, , \\
     \displaystyle\frac{\partial^3 u}{\partial n^3}  + (2 -\nu) \frac{\partial^3 u }{\partial n \partial \tau^2} + (1-\nu) \kappa \left( \frac{\partial^2 u}{\partial \tau^2} - \frac{\partial^2 u}{\partial n^2} \right)  = f_2 \, , &\qquad \text{ on } \partial \Omega \, ,
\end{cases} \label{freebcs}
\end{align}
where $\partial/\partial \tau$ denotes a tangential derivative and
$\kappa$ denotes the curvature of the boundary. As before, the first condition corresponds to a bending moment while the second condition corresponds to an equivalent shear force. These boundary conditions are 
common in applications where the motion at the edge of the plate is unrestricted
\cite{Ukrainskii2018,meylan96,meylan2002,sergienko13}.
}

There are many physical reasons to study exterior problems involving these boundary conditions. Plates that are found in buildings and bridges may be supported from within by internal columns, rods, and springs that affect their flexure \cite{zhao2002plate}. Frequently, cut-outs are made to these plates to lighten or ventilate the structure \cite{RAJAMANI1977549, HUANG1999769}, or to alter or eliminate resonant modes \cite{PhysRevB.73.064301, Lindsay2018}. This is valuable in the field of acoustics, where one is typically searching for a very specific frequency response \cite{ ZHU2024108814, laulagnet1998sound}. While the majority of our focus in this paper is on exterior problems, the techniques outlined here can also be applied for studying interior problems.

Traditionally, problems involving these boundary conditions have been treated using the finite element method (FEM) \cite{ monk87, climente2014gradient, mora2009piecewise, meylan2002}. For the {interior} clamped and supported plate problems, there are well-understood error estimates \cite{ishihara78, BRAMBLE1983} and software packages available for dealing with complex geometries. However, very few error estimates or high-order implementations exist for the free plate problem, which is typically solved on relatively simple geometries. {For exterior problems, implementation of na\"{i}ve} FEM {approaches face difficulties, due to the necessity of spatially truncating the problem. Typically, this involves incorporating auxiliary constraints at artificial boundaries. Two standard approaches for this are} transparent boundary conditions \cite{YUE2023100350, YUE2024112606} and perfectly matched layers \cite{FARHAT20112237}, though in either case it can be difficult to achieve zero reflectance for complicated wave phenomena. 

An alternative approach is to recast the PDE as an integral equation defined on the boundary of the domain. 
In general, this is possible whenever the PDE is linear and elliptic, with a known Green's function \cite{Atkinson_1997}. 
The basic approach is to define the solution of the PDE in terms of a layer potential with an unknown density
defined on the boundary of the domain. The layer potential is designed to automatically 
satisfy the PDE within the domain, while the unknown density is determined by enforcing the boundary conditions on the
boundary traces of the layer potential, resulting in a boundary integral equation (BIE). {When these boundary traces result in the sum of a bounded invertible operator and a compact operator, the BIE is referred to as a Fredholm equation of the second kind, or second kind integral equation (SKIE).} Numerically, BIEs offer several advantages over direct discretization or weak formulation of the
PDE. {In particular, for BIEs} all of the unknowns are located on the boundary of the domain, and hence only the boundary needs to be discretized. This is particularly useful for solving exterior problems, where any necessary 
radiation conditions can be incorporated into the definition of the layer potential{s}. {Additionally, it is well-known~\cite{kress}
that, under mild assumptions on the discretization scheme, the condition numbers of discrete matrix 
approximations of an invertible SKIE remain bounded as the boundary mesh is refined}. 

There is a significant literature for BIEs in the biharmonic {case} ($k=0$),
particularly for the clamped plate boundary conditions. SKIEs based on layer representations
were historically applied to analyze the existence, uniqueness, and regularity of solutions of the clamped plate 
problem~\cite{agmon1957multiple,cohen1983dirichlet}. An alternative route to the clamped plate problem
can be pursued by converting the problem to a Stokes problem, for which there are particularly
nice representations~\cite{muschelivsvili1933research,Rachh2018}.
A general approach for the $k=0$ case with both clamped plate and free plate problems
is developed in~\cite[\S 2.4]{hsiao2008boundary} {using suitable Calder\'{o}n 
projectors}, though the resulting boundary operators are not exactly SKIEs but include hypersingular {integral} operators. 
It is also possible to solve the supported plate problem by reduction
to a system of two Poisson equations  \cite{parisdeleon, sladek}, but this method requires volume integration and has no clear extension to the case when $k \neq 0$.

The thesis~\cite{farkas} describes a {general} framework for developing layer potential representations. The
basic idea is to {carefully} select linear combinations of derivatives of the Green's function {so} that the 
resulting BIE is second kind and {so that the singularities of the kernels cancel out when plugged into the highest order boundary condition}. In~\cite{farkas},
this technique is applied to the $k=0$ case for the clamped plate, fluid mechanics ($\nabla u$ is prescribed on the boundary),
and supported plate boundary conditions,
though with limited success in the case of the fluid mechanics and supported plate boundary conditions. 
In particular, the BIEs for the fluid mechanics and supported plate problems contain differential operators, {leading to poor conditioning upon discretization}. For the clamped plate problem, the framework of~\cite{farkas} extends readily
to three dimensions~\cite{JIANG20117488} and to nonzero $k$~\cite{Lindsay2018}.

In contrast with the steady-state ($k=0$) case, the oscillatory problem has received less attention, 
despite being a central focus of texts on structural engineering and elastic 
physics~\cite{timoshenko1959theory, landau59}. Previous attempts have relied on factoring the solution 
into Helmholtz and modified Helmholtz equations and solving a coupled system at the 
boundary~\cite{YUE2023100350, Dong2024, SMITH20114029}. 
These combined BIEs can suffer from ill-conditioning
though effective pre-conditioning strategies exist~\cite{Dong2024}. {Moreover}, these approaches are currently not capable of solving the free or supported plate problems, which motivates the {present work}. 

{The primary challenge in applying the approach of \cite{farkas} 
to the free plate problem lies in the order of the boundary conditions. Because
the derivatives in the boundary condition are high order, the integral kernels
in the representation can consist of at most first order derivatives of the Green's function.
In the case of the free plate, these kernels produce different types of singularities 
when plugged into the boundary conditions, so that it is difficult to obtain necessary
cancellations between terms.
In this paper, we present a method to circumvent this issue by composing standard layer 
potentials with the Hilbert transform, which has the effect of transforming the singularities 
into the necessary type. This approach offers 
new degrees of freedom in the design of the integral representations for high-order 
boundary conditions. The effect of the Hilbert transform can be understood in terms of 
the Poincar\'{e}-Bertrand formula~\cite{muskhelishvilisingularbook,shidong}, which
we use to establish that the resulting BIE is second kind.
A heuristic explanation of this technique is presented
in~\ref{app:heuristic}, using Fourier methods}.

The remainder of the paper is as follows. In Section {\ref{sec:mathframework}}, we give an overview of the Green's function and the boundary integral {equation} framework that will be used to solve this problem. In Section~\ref{clampedsection}, we review
the kernels previously used by~\cite{farkas, Lindsay2018} for solving the clamped plate problem and discuss their application to exterior {scattering} 
problems. In Section~\ref{supportedsection}, we give a representation for the supported plate which is similar to~\cite{farkas}, but uses extra terms to ensure that the BIE is second kind. In Section~\ref{freesection}, we present a novel representation for the free plate problem which uses the Hilbert transform in the representation of the solution. In Section~\ref{num_imp}, we discuss details of the implementation of these methods used for the
numerical results, describing a high-order-accurate approach to discretize the integral equations and stable methods to evaluate the kernels. In Section~\ref{results}, we present numerical results illustrating high-order convergence, and the application of these methods to solve scattering problems, including scattering from a large collection of objects. Finally, in Section~\ref{discussion}, we discuss the significance of the present work {and areas for future exploration}. 

\section{Mathematical framework}\label{sec:mathframework}

\label{biesection}

{
In this section, we outline our general framework for integral
representations of flexural wave boundary value problems.
We assume that the boundary of each subdomain, $\partial \Omega_i$, is a simple, regular curve 
with arc length parametrization $\bgamma(s)$, where $s \in [0,L]$ is the arc length parameter, 
and $L$ is the total arc length of the curve. The definitions and results presented here 
require that the curve is sufficiently regular, i.e. $\bgamma(s) \in C^{k}([0,L],\mathbb{R}^2)$ 
for $k$ sufficiently large. We also require that the curve is closed, so the 
parametrization satisfies the periodicity constraints $\bgamma^{(n)}(0) = \bgamma^{(n)}(L)$ for $n = 0, ..., k$. 
}

The Green's function for the free-space problem is defined as the solution to the equation
\begin{align}
    \Delta^2 u - k^4 u = \delta(\bx  - \by ) \, ,\quad \bx ,\by  \in  \mathbb{R}^2,
\end{align}
supplemented with the following radiation condition at infinity:
\begin{equation}
\label{eq:radcond}
   \frac{\bx}{\|\bx\|} \cdot \nabla u - i k u = o \left ( \frac{1}{\sqrt{\|\bx\|}} \right) \, . 
\end{equation}
{A straightforward calculcation shows that the Green's function for this problem is given by}
\begin{align}
    G(\bx ,\by ) &=  \frac{1}{2k^2} \left[ \frac{i}{4} H_0^{(1)} (k \|\bx  - \by \|) - \frac{1}{2\pi} K_0 (k \| \bx  - \by \|) \right] \, , \label{greens} 
\end{align}
 where $H_0^{(1)}$ is the zero-th order Hankel function of the first kind and $K_0$ is the zero-th order modified Bessel function of the second kind. 
 This Green's function is a scaled difference of the Green's functions for the Helmholtz and Yukawa equations, which oscillate and decay, respectively. As such, $G$ automatically satisfies the natural radiation condition \eqref{eq:radcond}. { Note that this Green's function has the same singularity near $\bx=\by$ as the biharmonic Green's function (see \ref{app:asymptotics}), which we denote by
\begin{align}
    G^B(\bx,\by) :=  \frac{1}{8\pi}\|\bx - \by\|^2 \ln \|\bx-\by \| \, .
\end{align} Therefore much of the analysis 
 of the integral operators will reduce to an analysis of $G^B$. }

In the BIE framework, we seek to represent the solution to a boundary value problem (BVP) as a combination of
{layer potentials} applied to some densities on the boundary:
\begin{align}
\label{eq:layerpotential}
    u(\bx ) &= \mathcal{K}_1 [\rho_1 ](\bx ) + \mathcal{K}_2 [\rho_2 ](\bx ) \, ,
\end{align}
where $\rho_1 $ and $\rho_2 $ are densities that lie in a suitable function space on the 
boundary (e.g. $L^2(\partial \Omega)$) and where the integral operators $\mathcal{K}_1$ and $\mathcal{K}_2$ are given by
\begin{align}
    \mathcal{K}_1[\rho_1](\bx ) &:= \int_{\partial \Omega} K_1(\bx ,\by ) \rho_1(\by )  \, \dd S(\by ) \, ,  \\
    \mathcal{K}_2[\rho_2](\bx ) &:= \int_{\partial \Omega} K_2(\bx ,\by ) \rho_2(\by )  \, \dd S(\by )  \, .
\end{align}

In this formalism, we represent the solution of the BVP using a collection of ``charges'' (sources, dipoles, and multipoles) defined on the boundary $\partial\Omega$. Since the kernels $K_1$ and $K_2$ {involve carefully chosen linear
combinations of the Green’s function and its derivatives,} the representation automatically satisfies \eqref{flexural} 
{for points off of the boundary. What remains is to ensure consistency with the boundary conditions. In order to obtain an equation for the densities, we insert our representation into the boundary conditions by taking a limit of the appropriate linear combination of derivatives of our ansatz $u$ as the target approaches the boundary from a normal direction.} This leads to a $2\times 2$ system of integral  equations on the boundary $\partial \Omega$:
\begin{equation}
\label{eq:biesystem}
    \begin{pmatrix}
    D_{11}{(\bx)} & D_{12}{(\bx)} \\
    D_{21}{(\bx)} & D_{22}{(\bx)}
    \end{pmatrix}
    \begin{pmatrix}
        \rho_1(\bx ) \\
        \rho_2(\bx ) 
    \end{pmatrix} + \pv \int_{\partial \Omega} \begin{pmatrix}
    K_{11}(\bx ,\by ) & K_{12}(\bx ,\by ) \\
    K_{21}(\bx ,\by ) & K_{22}(\bx ,\by )
    \end{pmatrix}
    \begin{pmatrix}
        \rho_1(\by ) \\
        \rho_2(\by )
    \end{pmatrix} \, \dd S(\by )
    = \begin{pmatrix}
        f_1(\bx ) \\ f_2(\bx )
    \end{pmatrix} \, ,
\end{equation}
where {the expressions $f_1$ and $f_2$ on the right-hand side depend on the boundary data} and $\pv$ indicates that the (possibly singular) integral is to be understood in the principal
value sense. It is useful to define the on-surface integral operators as well:
\begin{align}
    \mathcal{K}_{ij}[\rho_j]({\bx_0} ) &:= \pv \int_{\partial \Omega} K_{ij}({\bx_0} ,\by ) \rho_j(\by )  \, \dd S(\by ) {\, , \quad \bx_0 \in \partial \Omega} \, .
\end{align}
The jump {values, $D_{ij}$,} 
are precisely the difference between the limit of the integral operator as it approaches the boundary and the operator on the boundary:
    \begin{align}
        D_{ij}{(\bx_0)} \rho_j(\bx _0) :=  \lim_{\bx  \to  \bx _0} \int_{\partial\Omega} K_{ij}(\bx ,\by ) \rho_j(\by ) \, \dd S(\by ) - \pv \int_{\partial\Omega} K_{ij}(\bx _0,\by ) \rho_j(\by ) \, \dd S(\by ) \, .
    \end{align}
    {When the system \eqref{eq:biesystem} is in the form of a bounded and invertible
    operator plus a compact operator}, the equation is referred to as a {Fredholm equation of the second kind,
    or second kind integral equation (SKIE)}. 
    The {aim of the present work is to select} kernels $K_1$ and $K_2$ for each of the clamped plate, 
    free plate, and supported
    plate boundary value problems so that system \eqref{eq:biesystem} is an SKIE. 
    {For details on heuristic methods for deriving such kernels, see~\ref{app:heuristic}.}

{
    \begin{rmk}
        For the sake of simplicity, the jump matrix is taken to be a multiplication
        operator above. As we will see in Section~\ref{freesection}, the jump
        matrix can in general contain surface operators. Further, we need
        to consider Hadamard finite part integrals, in addition to Cauchy 
        principal value integrals, in some cases.
    \end{rmk}}
    
In the remainder of this manuscript, we use the following conventions for notation. For a point $\bx$ on $\partial \Omega$, $\bn (\bx)$ is the {normal that points into the exterior},
$\btau (\bx)$ is the positively-oriented (clockwise) unit tangent vector, $\kappa(\bx)$ is 
the signed curvature at the point $\bx$, {and $\kappa'(\bx)$ is 
the arc length derivative of curvature. For the sake of brevity, we will
often write, in a standard abuse of notation, $\kappa(s), \bn(s)$, and $\btau(s)$ in place
of $\kappa(\bgamma(s)), \bn(\bgamma(s))$, and  $\btau(\bgamma(s))$, respectively, while performing Taylor expansions of the kernels.} 
The input to the Green's function is usually of the form $G(\bx,\by)$, where $\bx$ is 
sometimes referred to as the ``target'' and $\by$ as the ``source.'' We often employ the 
more compact notation $\bn_\bx  = \bn (\bx)$, $\btau_\by =\btau(\by)$, etc. We denote directional derivatives of the Green's function by subscripts, e.g. 
$G_{\bn_\bx}(\bx,\by) = \bn (\bx) \cdot \nabla_\bx  G(\bx,\by)$ and $G_{\btau_\by} = \btau(\by) \cdot \nabla_\by G$. 
Multiple subscripts correspond to contractions of the Hessian or other higher-order
derivative tensors of $G$ with the indicated directions, e.g. 
\begin{equation}
    G_{\bn_\bx \btau_\bx \btau_\by}(\bx,\by) = \bn_\bx \cdot \nabla_\mathbf{w} \left ( \btau_\bx \cdot \nabla_\mathbf{w} \left ( \btau_\by \cdot \nabla_\mathbf{z} G(\mathbf{w},\mathbf{z})
    \right ) \right )  |_{\mathbf{w}=\bx,\mathbf{z}=\by} \, .
\end{equation}
{Writing $\bx(s) := \bgamma(s)$, } an important distinction is that
\begin{equation}
\frac{\dd}{\dd s} G(\bx(s),\by) = G_{\btau_\bx}(\bx(s),\by) \, , \quad \textrm{but} \quad 
\frac{\dd^2}{\dd s^2} G(\bx(s),\by) = G_{\btau_\bx \btau_\bx}(\bx(s),\by) -
\kappa(\bx(s)) G_{\bn_\bx}(\bx(s),\by) \, ,
\end{equation}
which can be easily seen through application of the product rule. 

\section{The clamped plate problem} \label{clampedsection}

As noted in the introduction, \cite{farkas} developed kernels for the biharmonic equation ($k = 0$)
by selecting linear combinations of derivatives of the Green's function to optimize the regularity
of the resulting boundary integral equation kernels. These same linear combinations apply readily to
the case of non-zero $k$ \cite{Lindsay2018}. 
For the sake of completeness, we {restate} these kernels here:
\begin{align}
    K_1 &= G_{\bn_\by  \bn_\by  \bn_\by } + 3G_{\bn_\by  \btau_\by  \btau_\by } \label{clampedk1} , \\
    K_2 &= - G_{\bn_\by  \bn_\by } + G_{\btau_\by  \btau_\by } \, . \label{clampedk2}
\end{align}
Below are the kernels that appear in the boundary integral equation:
\begin{align}
    K_{11} &= G_{\bn_\by  \bn_\by  \bn_\by } + 3G_{\bn_\by  \btau_\by  \btau_\by } \label{first_clamped_ker}  , \\
    K_{12} &= - G_{\bn_\by  \bn_\by } + G_{\btau_\by  \btau_\by }  , \\
    K_{21} &= G_{\bn_\bx  \bn_\by  \bn_\by  \bn_\by } + 3G_{\bn_\bx  \bn_\by  \btau_\by  \btau_\by } , \\
    K_{22} &= - G_{\bn_\bx  \bn_\by  \bn_\by } + G_{\bn_\bx  \btau_\by  \btau_\by }  . \label{last_clamped_ker}
\end{align}
{
\begin{thm}
\label{thm:clamped}
    Let $\cK_1$ and $\cK_2$ be the layer potentials corresponding to the
    kernels $K_1$ and $K_2$, respectively. Let $\cK_{11},\cK_{12},\cK_{21},$ and 
    $\cK_{22}$ 
    be the boundary operators corresponding to the kernels $K_{11},K_{12},K_{21}$, 
    and $K_{22}$, respectively. These operators satisfy the following jump relations:
\begin{align}
    \lim_{\bx \to \bx_0^{\pm}}\mathcal{K}_{1}[\sigma](\bx) &= \mp \frac{1}{2} \sigma(\bx_0) + \mathcal{K}_{11}[\sigma](\bx_0) \, , \\
    \lim_{\bx \to \bx_0^{\pm}}\mathcal{K}_{2}[\sigma](\bx) &= \mathcal{K}_{12}[\sigma](\bx_0) \, , \\
    \lim_{\bx \to \bx_0^{\pm}}\bn(\bx_0) \cdot \nabla_{\bx}  \mathcal{K}_{1}[\sigma](\bx) &=  \pm \kappa(\bx_0) \sigma(\bx_0) + \mathcal{K}_{21}[\sigma](\bx_0) \, , \\
    \lim_{\bx \to \bx_0^{\pm}}\bn(\bx_0) \cdot \nabla_{\bx} \mathcal{K}_{2}[\sigma](\bx) &= \mp \frac{1}{2} \sigma(\bx_0) + \mathcal{K}_{22}[\sigma](\bx_0) \, ,
\end{align}
for $\bx_0 \in \partial\Omega$, where $\bx_0^+$ corresponds to the limit from the exterior and $\bx_0^-$ corresponds to the limit from the interior, each along a normal direction.

Moreover, the corresponding exterior and interior integral equations
for the representation $u = \cK_1[\rho_1] + \cK_2[\rho_2]$, i.e. 
\begin{equation}
    \begin{pmatrix}
        \mp \frac{1}{2} I & 0\\
        \pm \kappa & \mp\frac{1}{2} I
    \end{pmatrix} \begin{pmatrix} \rho_1 \\ \rho_2 \end{pmatrix} + \begin{pmatrix} 
    \cK_{11} & \cK_{12} \\ 
    \cK_{21} & \cK_{22}
    \end{pmatrix} \begin{pmatrix} \rho_1 \\ \rho_2 \end{pmatrix} = \begin{pmatrix} f_1 \\ f_2 \end{pmatrix} \, ,
\end{equation}
are second kind on $L^2(\partial \Omega) \times L^2(\partial \Omega)$. 
\end{thm}
}

{A proof of the above theorem is contained in \cite{farkas} for the $k=0$ case
and the arguments extend readily to the $k\ne 0$ case. We reproduce similar arguments here 
so that the exposition is self-contained. Before proving the theorem, we establish the 
continuity of the kernels $K_{ij}$ for sources and targets on the boundary.}

\begin{lem}
The kernel functions $K_{11}$, $K_{12}$, $K_{21}$, and $K_{22}$ are continuous when
restricted to the boundary. \label{lem:clampedsmooth}
\end{lem}
{
\begin{proof}
    We prove this for the case that $\Omega$ is a single, simply connected domain;
    the extension to the case of multiple components is straightforward.
    Based on the asymptotic expansion of the Green's function (\ref{app:asymptotics}), it is sufficient 
    for us to show the continuity of the corresponding kernels with $G$ replaced by
    the biharmonic Green's function $G^B$,
    which is the dominant part of $G$ as $\|\bx -\by \| \to 0$. A similar observation was made 
    in \cite{Lindsay2018} for a variant of the flexural Green's function considered here. 
    
    We denote these kernels by $K_{ij}^B$. Because $G$ is smooth for $\bx \ne \by$,
    it is sufficient to show continuity as $\|\bx-\by\|\to 0$ for $\bx,\by \in \partial \Omega$.
    First, we perform a Taylor expansion of the 
    arc length parametrization $\bgamma(t+s)$ about some $\bgamma(t)$: 
     \begin{align}
          \bgamma(t + s) &= \bgamma(t) + s \btau(t) - \frac{s^2}{2} \kappa(t) \bn(t)  - \frac{s^3}{6}(\kappa'(t) \bn(t) + \kappa^2(t) \btau(t)) + \mathcal{O}(s^4) \, , \label{gammaseries}
     \end{align}
    where $\kappa(t)$ and $\kappa'(t)$ are the curvature and arc length derivative of curvature at $\bgamma(t)$. Because integration is performed with respect to $\by$ for fixed $\bx$, we let $ \by  = \bgamma(s) $ and fix $\bx  = \bgamma(0) $ for convenience. Setting $t = 0$ in the formula above: 
     \begin{align}
          \by  &= \bx + s\btau(0) - \frac{s^2}{2} \kappa(0) \bn(0)  - \frac{s^3}{6}(\kappa'(0) \bn(0) + \kappa^2(0) \btau(0))  + \mathcal{O}(s^4) \, . \label{yseries}
     \end{align}
    Moreover, we can expand $  \btau(s) $ and $ \bn(s) $: 
    \begin{align}
         \btau (s)  &= \btau(0) - s \kappa(0) \bn(0) -\frac{s^2}{2}(\kappa'(0) \bn(0) + \kappa^2(0) \btau(0))  + \mathcal{O}(s^3) \, , \label{tauseries}\\
         \bn (s)  &= \bn(0) + s \kappa(0) \btau(0)   + \frac{s^2}{2}(\kappa'(0) \btau(0) - \kappa^2(0) \bn(0))  + \mathcal{O}(s^3) \, . \label{nseries}
    \end{align}
    We first analyze the asymptotics of the kernel $K_{11}^{B} :=  G^{B}_{\bn_\by  \bn_\by  \bn_\by } + 3 G^{B}_{\bn_\by  \btau_\by  \btau_\by } $. After inserting the above expansions into the formulae derived in \ref{appclamped}, we obtain:
    \begin{align}
        G^{B}_{\bn_\by  \bn_\by  \bn_\by }(\bgamma(0),\bgamma(s)) &= \frac{3}{8\pi}\kappa(0) + \mathcal{O}(s) \, , \\
        G^{B}_{\bn_\by  \btau_\by  \btau_\by } (\bgamma(0),\bgamma(s)) &= -\frac{1}{8\pi}\kappa(0) + \mathcal{O}(s) \, .
    \end{align}
    Therefore, the on-surface limit of $K^B_{11}(\bx,\by)$ as $\by$ approaches $\bx$ along the boundary is
    \begin{align}
        \lim_{s\to  0} K^B_{11}(\bgamma(t),\bgamma(t+s)) = 0 \, .
    \end{align}
    We conclude that $K^B_{11}$ is continuous. Next, we examine $K^B_{12} := -G^B_{\bn_\by \bn_\by} + G^B_{\btau_\by \btau_\by}$. Using the formulae in \ref{appclamped} we have:
    \begin{align}
        K^B_{12}(\bx,\by) = \frac{1}{4\pi} \frac{[(\bx - \by)\cdot  \bn (\by) ]^2}{\|\bx-\by\|^2} - \frac{1}{4\pi} \frac{[(\bx-\by)\cdot \btau (\by)]^2}{\|\bx-\by\|^2} \, .
    \end{align}
    Expanding this kernel as above, we find that:
    \begin{align}
        \lim_{s\to  0} K^B_{12}(\bgamma(t),\bgamma(t+s)) = \frac{1}{4\pi} \, .
    \end{align}
    Therefore $K^B_{12}$ is continuous as well. We now examine the kernels that arise from applying the second boundary condition, i.e. $K^B_{21}$ and $K^B_{22}$. Looking at the terms in $K^B_{22} := - G^B_{\bn_\bx \bn_\by \bn_\by } + G^B_{\bn_\bx \btau_\by \btau_\by }$, we have:
    \begin{align}
        G^B_{\bn_\bx \bn_\by \bn_\by }(\bgamma(0),\bgamma(s)) &= - \frac{1}{8\pi} \kappa(0) + \mathcal{O}(s) \, , \\
        G^B_{\bn_\bx \btau_\by \btau_\by }(\bgamma(0),\bgamma(s)) &= \frac{3}{8\pi} \kappa(0) + \mathcal{O}(s) \, .
    \end{align}
    The on-surface limit for $K^B_{22}$ becomes:
    \begin{align}
        \lim_{s\to  0} K^B_{22}(\bgamma(t),\bgamma(t+s)) = \frac{1}{2\pi}\kappa(t) \, .
    \end{align}
    Lastly, expanding the two kernels in $K_{21}^B := G^B_{\bn_\bx  \bn_\by  \bn_\by  \bn_\by } + 3G^B_{\bn_\bx  \bn_\by  \btau_\by  \btau_\by }$, we get:
    \begin{align}
        G^B_{\bn_\bx  \bn_\by  \bn_\by  \bn_\by }(\bgamma(0),\bgamma(s)) &= - \frac{3}{4\pi}\frac{1}{s^2} + \frac{5}{16 \pi} \kappa^2(0)  + \mathcal{O}(s) \, , \\
        G^B_{\bn_\bx  \bn_\by  \btau_\by  \btau_\by}(\bgamma(0),\bgamma(s)) &= \frac{1}{4\pi}\frac{1}{s^2} - \frac{17}{48\pi} \kappa^2(0) + \mathcal{O}(s) \, .
    \end{align}
    Therefore, the on-surface limit of $K_{21}^B$ becomes:
    \begin{align}
        \lim_{s\to  0} K_{21}^B(\bgamma(t),\bgamma(t+s)) = -\frac{3}{4\pi} \kappa^2(t) \, .
    \end{align} 
    Provided that $\partial \Omega$ is sufficiently smooth, $\kappa(t)$ is continuous, therefore the kernels $K_{21}^B$ and $K_{22}^B$ are continuous on $\partial \Omega$. The continuity of $K_{11}, K_{12}, K_{21},$ and $K_{22}$ follows immediately. 
\end{proof} 
}

{\begin{rmk}
    We observe that even though $G_{\bn_\bx  \bn_\by  \bn_\by  \bn_\by }$ and $G_{\bn_\bx  \bn_\by  \btau_\by  \btau_\by}$ are both singular, by a judicious choice of coefficients in $K_{1},$ these singularities cancel in $K_{21}.$ The kernels for solving both the free and supported plate problems will similarly rely on cancellations of this type.
\end{rmk}}

\begin{proof}[Proof of Theorem~\ref{thm:clamped}.]
    {The jump properties are established in~\ref{app:jumpcondclamped}.}
        {It is sufficient for us to show that the integral operators $\mathcal{K}_{11}, \mathcal{K}_{12}, \mathcal{K}_{21},$ and $ \mathcal{K}_{22}$ are compact.  It is known from the theory of integral equations \cite{kress} that integral operators on $\partial \Omega$ with continuous kernels are compact from $L^2(\partial \Omega) \to L^2(\partial \Omega)$. 
        The result then follows from Lemma~\ref{lem:clampedsmooth}.}
\end{proof}

\section{The supported plate problem} \label{supportedsection}

For {the supported plate problem  \eqref{supportedbcs},} we use a representation 
that is similar to \cite{farkas} with several additional terms that ensure 
{the resulting boundary integral equation is second kind}:
\begin{align}
    K_{1} &= G_{\bn_\by  \bn_\by  \bn_\by }   + \alpha_1 G_{\bn_\by  \btau_\by  \btau_\by }  +  \alpha_2 \kappa(\by ) G_{\bn_\by  \bn_\by  } + \alpha_3 \kappa'(\by ) G_{\btau_\by  } \, ,  \\
    K_{2} &= G_{\bn_{\by }} \, ,
\end{align}
where {$\kappa'(\by)$ is the arc length derivative of curvature at $\by \in \partial \Omega$} and the coefficients $\alpha_1, \alpha_2,$ and $ \alpha_3$ are given by:
\begin{align}
    \alpha_1 &= 2-\nu \, , \\
    \alpha_2 &= \frac{(-1+\nu)(7+\nu)}{3-\nu} \, , \\
    \alpha_3 &= \frac{(1-\nu)(3+\nu)}{1+\nu}  \, .
\end{align}

The corresponding kernel functions in the integral equation are:
\begin{align}
    K_{11} &= G_{\bn_\by  \bn_\by  \bn_\by }   + \alpha_1 G_{\bn_\by  \btau_\by  \btau_\by }  +  \alpha_2 \kappa(\by ) G_{\bn_\by  \bn_\by  } + \alpha_3 \kappa'(\by ) G_{\btau_\by  } \, , \\
    K_{12} &= G_{\bn_{\by }} \, , \\
    K_{21} &=  G_{\bn_\bx  \bn_\bx  \bn_\by  \bn_\by  \bn_\by } + \nu  G_{\btau_\bx  \btau_\bx  \bn_\by  \bn_\by  \bn_\by }   + \alpha_1 G_{ \bn_\bx  \bn_\bx  \bn_\by  \btau_\by  \btau_\by } + \nu  \alpha_1 G_{\btau_\bx  \btau_\bx  \bn_\by  \btau_\by  \btau_\by } + \nonumber  \\
    & \qquad  \alpha_2 \kappa(\by ) G_{\bn_\bx  \bn_\bx  \bn_\by  \bn_\by  } +\nu \alpha_2 \kappa(\by ) G_{\btau_\bx  \btau_\bx  \bn_\by  \bn_\by  } + \alpha_3 \kappa'(\by ) G_{\bn_\bx  \bn_\bx  \btau_\by  } + \nu \alpha_3 \kappa'(\by ) G_{\btau_\bx  \btau_\bx  \btau_\by  } \, ,  \\
    K_{22} &= G_{\bn_\bx  \bn_\bx  \bn_{\by }}  + \nu  G_{\btau_\bx  \btau_\bx  \bn_{\by }}  \, .\label{last_supported_kernel}
\end{align}

{ 
\begin{thm}
\label{thm:supported}
    Let $\cK_1$ and $\cK_2$ be the layer potentials corresponding to the
    kernels $K_1$ and $K_2$, respectively. Let $\cK_{11},\cK_{12},\cK_{21},$ and 
    $\cK_{22}$ 
    be the boundary operators corresponding to the kernels $K_{11},K_{12},K_{21}$, 
    and $K_{22}$, respectively. Let 
    $$ c_0 = \frac{(\nu -1) (\nu +3) (2 \nu -1)}{2 (3-\nu)} \, .$$ 
    These operators satisfy the following jump relations:
\begin{align}
    \lim_{\bx \to \bx_0^{\pm}}\mathcal{K}_{1}[\sigma](\bx) &= \mp \frac{1}{2} \sigma(\bx_0) + \mathcal{K}_{11}[\sigma](\bx_0) \, , \\
    \lim_{\bx \to \bx_0^{\pm}}\mathcal{K}_{2}[\sigma](\bx) &= \mathcal{K}_{12}[\sigma](\bx_0) \, , \\
    \lim_{\bx \to \bx_0^{\pm}}\left ( \nu \Delta + (1-\nu) \left (\bn(\bx_0) \cdot \nabla_{\bx} \right )^2
    \right ) \mathcal{K}_{1}[\sigma](\bx) &=  \pm c_0 \kappa^2(\bx_0) \sigma(\bx_0) + \mathcal{K}_{21}[\sigma](\bx_0) \, , \\
    \lim_{\bx \to \bx_0^{\pm}}\left ( \nu \Delta + (1-\nu) \left ( \bn(\bx_0) \cdot \nabla_{\bx}\right)^2
    \right ) \mathcal{K}_{2}[\sigma](\bx) &= \mp \frac{1}{2} \sigma(\bx_0) + \mathcal{K}_{22}[\sigma](\bx_0) \, ,
\end{align}
for $\bx_0 \in \partial\Omega$, where $\bx_0^+$ corresponds to the limit from the exterior and $\bx_0^-$ corresponds to the limit from the interior, each along a normal direction.

Moreover, the corresponding exterior and interior integral equations for the 
representation $u = \cK_1[\rho_1] + \cK_2[\rho_2]$, i.e. 
\begin{equation}
    \begin{pmatrix}
        \mp \frac{1}{2} I & 0\\
        \pm c_0 \kappa^2 & \mp\frac{1}{2} I
    \end{pmatrix} \begin{pmatrix} \rho_1 \\ \rho_2 \end{pmatrix} + \begin{pmatrix} 
    \cK_{11} & \cK_{12} \\ 
    \cK_{21} & \cK_{22}
    \end{pmatrix} \begin{pmatrix} \rho_1 \\ \rho_2 \end{pmatrix} = \begin{pmatrix} f_1 \\ f_2 \end{pmatrix} \, ,
\end{equation}
are second kind on $L^2(\partial \Omega) \times L^2(\partial \Omega)$
for $\nu \not\in \{-1,3\}$. 
\end{thm}}

{Before proving the above theorem, we require the following lemma.}

\begin{lem}
For $\nu \notin \{-1, 3\}$, the kernels $K_{11},K_{12},K_{21},$ and $K_{22}$  are continuous. \label{supportedsmooth}
\end{lem}
\begin{proof}
    We prove this for the case that $\Omega$ is a single, simply connected domain;
    the extension to the case of multiple components is straightforward.
As in the proof of Lemma~\ref{lem:clampedsmooth}, let $G^{B}$ be the biharmonic Green's
function and let $K^B_{ij}$ be the integral kernels above with $G$ replaced by 
$G^B$. Since these kernels involve at most five derivatives and $G$ is smooth for 
$\bx \ne \by$, {it is} sufficient for us 
to show the continuity of $ K^{B}_{11}, K^B_{12}, K^B_{21}, K^B_{22}$ as $\|\bx -\by \| \to 0$
for $\bx,\by \in \partial \Omega$. 
{ Again, we let $\bx = \bgamma(0)$ and $\by = \bgamma(s)$, and we let $\kappa(t)$ and $\kappa'(t)$ be the curvature and arc length derivative of curvature at point $\bgamma(t)$. Because our kernels now include five derivatives, we require one extra term in the expansions of $\by,\bn,$ and $\btau$, i.e.
\begin{align}
          \by  &= \bx + s\btau(0) - \frac{s^2}{2} \kappa(0) \bn(0)  - \frac{s^3}{6}(\kappa'(0) \bn(0) + \kappa^2(0) \btau(0)) \nonumber \\
          &\qquad \qquad\qquad\qquad\qquad\qquad + \frac{s^4}{24}((- \kappa''(0) + \kappa^3(0)) \bn(0) - 3\kappa(0) \kappa'(0) \btau(0) ) + \mathcal{O}(s^5) \, , \label{expansion1} \\
         \btau (s)  &= \btau(0) - s \kappa(0) \bn(0) -\frac{s^2}{2}(\kappa'(0) \bn(0) + \kappa^2(0) \btau(0)) \nonumber \\
         &\qquad\qquad\qquad\qquad\qquad  +\frac{s^3}{6}((- \kappa''(0) + \kappa^3(0)) \bn(0) - 3\kappa(0) \kappa'(0) \btau(0) ) + \mathcal{O}(s^4) \, , \\
         \bn (s)  &= \bn(0) + s \kappa(0) \btau(0)   + \frac{s^2}{2}(\kappa'(0) \btau(0) - \kappa^2(0) \bn(0)) \nonumber \\
         &\qquad\qquad\qquad\qquad\qquad + \frac{s^3}{6} ((\kappa''(0) - \kappa^3(0))\btau(0)-3\kappa(0) \kappa'(0) \bn(0)) + \mathcal{O}(s^4) \, . 
    \end{align} }

The kernels in $K^B_{11}, K^B_{12},$ and $K^B_{22}$ are similar
to those in Lemma~\ref{lem:clampedsmooth}, so we do not repeat the analogous arguments here. 
It remains to show that $K^B_{21}$ is continuous. The curvature of the boundary at the source point can be expanded as:
\begin{align}
    \kappa( s) = \kappa(0) + s \kappa'(0) + \frac{1}{2} s^2 \kappa''(0) +  \mathcal{O}(s^3) \, ,
\end{align}
while the arc length derivative of curvature can also be expanded: 
\begin{align}
    \kappa'( s) = \kappa'(0) + s \kappa''(0) +  \mathcal{O}(s^2) \, . \label{expansion5}
\end{align}
Using {expansions \eqref{expansion1}-\eqref{expansion5} together with the formulae for derivatives of the Green's function in~\ref{appsupported}}, we obtain the on-surface asymptotics for the terms in $K^B_{21}$:
\begin{align}
        G^{B}_{\bn_\bx  \bn_\bx  \bn_\by  \bn_\by  \bn_\by  } {(\bgamma(0),\bgamma(s))} &= - \frac{3 \kappa(0)}{4\pi} \frac{1}{s^2} - \frac{ \kappa'(0)}{ \pi } \frac{1}{s} + \frac{3 \kappa^3(0) - 7 \kappa''(0)}{16 \pi} + \mathcal{O}(s) \, , \\
        G^{B}_{\btau_\bx  \btau_\bx  \bn_\by  \bn_\by  \bn_\by  } {(\bgamma(0),\bgamma(s))} &= -\frac{3 \kappa(0)}{4\pi} \frac{1}{s^2} + \frac{3 \kappa^3(0) + \kappa''(0)}{16 \pi} + \mathcal{O}(s) \, , \\
        G^B_{ \bn_{\bx } \bn_{\bx } \bn_\by  \btau_\by  \btau_\by } {(\bgamma(0),\bgamma(s))} &= \frac{5\kappa(0)}{4\pi}\frac{1}{s^2} + \frac{\kappa'(0)}{\pi} \frac{1}{s}+ \frac{-7\kappa^3(0) + 19\kappa''(0)}{48\pi} + \mathcal{O}(s) \, , \\
        G^B_{\btau_\bx  \btau_\bx  \bn_\by  \btau_\by  \btau_\by } {(\bgamma(0),\bgamma(s))} &= \frac{\kappa(0)}{4\pi} \frac{1}{s^2} - \frac{11 \kappa^3(0) + \kappa''(0)}{48 \pi} + \mathcal{O}(s) \, , \\
        \kappa( {s} ) \, G^B_{\bn_\bx  \bn_\bx  \bn_\by  \bn_\by } {(\bgamma(0),\bgamma(s))} &= \frac{3 \kappa(0)}{4 \pi  } \frac{1}{s^2} +\frac{3 \kappa'(0)}{4 \pi  }\frac{1}{s} +\frac{6 \kappa''(0) - 3 \kappa^3(0)}{16 \pi } + \mathcal{O}(s) \, , \\
        \kappa( {s}  ) \, G^B_{\btau_\bx  \btau_\bx  \bn_\by  \bn_\by } {(\bgamma(0),\bgamma(s))} &=  -\frac{\kappa(0)}{4 \pi  }\frac{1}{s^2}-\frac{\kappa'(0)}{4 \pi }\frac{1}{s} -\frac{6\kappa''(0) - \kappa^3(0)}{48 \pi } + \mathcal{O}(s) \, , \\
        \kappa'( {s} ) \, G^B_{\bn_{\bx } \bn_\bx  \btau_\by } {(\bgamma(0),\bgamma(s))} &=  \frac{\kappa'(0)}{4 \pi  }\frac{1}{s} + \frac{\kappa''(0)}{4 \pi } + \mathcal{O}(s)  \, , \\
        \kappa'( {s}  ) \, G^B_{\btau_{\bx } \btau_\bx  \btau_\by } {(\bgamma(0),\bgamma(s))} &=  \frac{\kappa'(0)}{4 \pi  }\frac{1}{s} + \frac{\kappa''(0)}{4 \pi } + \mathcal{O}(s) \, .
    \end{align} 
It follows that:
\begin{align}
    \lim_{s \to  0} K^B_{21} {(\bgamma(t) ,\bgamma(t+s) )} &=  \frac{(\nu-1) \left(12 \kappa^3(t) \left(\nu^2-\nu+4\right)+\kappa''(t) \left(-5 \nu^2+4 \nu+33\right)\right)}{48 \pi  (\nu-3)} \, .
\end{align}
{Provided that the boundary $\partial \Omega$ is sufficiently smooth,} this limit is well-defined, and the kernel $K^B_{21}$ is continuous. 
\end{proof}

{
\begin{proof}[Proof of Theorem~\ref{thm:supported}.]
The jump properties are established in \ref{app:jumpcondsupp}. It is sufficient
then to show that the operators $\cK_{11}$, $\cK_{12}$, $\cK_{21}$,
and $\cK_{22}$ are compact. 
By the same reasoning as in the proof of Theorem~\ref{thm:clamped}, 
compactness follows from Lemma~\ref{supportedsmooth}.
\end{proof}
}

\section{The free plate problem} \label{freesection}

Before proceeding, we recall the definitions of {two common} 
integral operators: the Hilbert transform and {the} Laplace double layer, {defined}
for a {closed} curve $\partial \Omega$.
The Hilbert transform $\mathcal{H}$ of a function $f$ on $\partial \Omega$ is given by:
\begin{align}
    \mathcal{H}[f](\bx ) &=  \int_{\partial \Omega}  \frac{(\bx -\by ) \cdot  \btau (\by) }{ {\pi} \|\bx -\by \|^2} f(\by ) \, \dd S(\by ) \, , \label{hilbertdef}
\end{align}
{while} the Laplace double layer potential $\mathcal{D}$ of a function $f$ on $\partial \Omega$  is given by:
\begin{align}
    \mathcal{D}[f](\bx ) &=  \int_{\partial \Omega}  \frac{(\bx -\by ) \cdot  \bn (\by) }{{2\pi}\|\bx -\by \|^2} f(\by ) \, \dd S(\by )  \, . \label{doublelayer}
\end{align}

{When $\Omega$ is given as a single, simply-connected domain, we represent} 
the solution {of the free plate problem} \eqref{freebcs} using 
the layer potential operators:
\begin{align}
    \mathcal{K}_1 &= \mathcal{K}_1^a + {\beta^\pm} \mathcal{K}_1^b  \mathcal{H} \, , \\
    \mathcal{K}_2 &= \mathcal{K}_2 \, ,
\end{align}
where $\beta := \displaystyle\frac{1+\nu}{2}$ { and $\beta^\pm := \pm \beta$. The coefficient $\beta^+$ is used for the exterior problem and $\beta^-$ is used for the interior problem.} The kernels of the integral operators above are given by:
\begin{align}
    K^a_{1}  &= G_{\bn_\by } \, , \\
    K^b_{1} &=  G_{\btau_\by } \, , \\
    K_{2} &= G \, .
\end{align}
{Taking suitable derivatives results in the following boundary integral operators:
\begin{align}
    \mathcal{K}_{11} &= \mathcal{K}_{11}^a + {\beta^\pm} \mathcal{K}_{11}^b  \mathcal{H} \, , \\
    \mathcal{K}_{12} &= \mathcal{K}_{12} \, , \\
    \mathcal{K}_{21} &= \mathcal{K}_{21}^a + {\beta^\pm} \mathcal{K}_{21}^b  \mathcal{H} \, , \\
    \mathcal{K}_{22} &= \mathcal{K}_{22} \, ,
\end{align}
where,
\begin{align}
    K^a_{11} &= G_{\bn_\bx \bn_\bx \bn_\by } + \nu G_{\btau_\bx \btau_\bx \bn_\by } \, , \label{first_free_ker} \\
    K^b_{11} &=  G_{\bn_\bx \bn_\bx \btau_\by } + \nu  G_{\btau_\bx \btau_\bx \btau_\by } \, , \\
    K_{12} &= G_{\bn_\bx \bn_\bx } + \nu G_{\btau_\bx \btau_\bx } \, , \\
    K_{21}^a &=  G_{\bn_\bx  \bn_\bx  \bn_\bx \bn_\by} + (2-\nu) G_{\bn_\bx  \btau_\bx  \btau_\bx  \bn_\by} + (1-\nu)\kappa(\bx) \left(G_{\btau_\bx  \btau_\bx \bn_\by} -  G_{\bn_\bx  \bn_\bx \bn_\by}\right) \, , \\
    K_{21}^b &=  G_{\bn_\bx  \bn_\bx  \bn_\bx \btau_\by} +  (2-\nu) G_{\bn_\bx  \btau_\bx  \btau_\bx \btau_\by } +  (1-\nu)\kappa(\bx) \left(G_{\btau_\bx  \btau_\bx \btau_\by} -  G_{\bn_\bx  \bn_\bx \btau_\by}\right) \, , \\
    K_{22} &= G_{\bn_\bx  \bn_\bx  \bn_\bx } + (2-\nu) G_{\bn_\bx  \btau_\bx  \btau_\bx  } + (1-\nu)\kappa(\bx) \left(G_{\btau_\bx  \btau_\bx } -  G_{\bn_\bx  \bn_\bx }\right)  \, .
\end{align} 
}
{
\begin{rmk}
    Practitioners of integral equations may be surprised that 
    the interior and exterior representations exhibit such a 
    sign dependence. The reason is that there is a cancellation
    which occurs between the $\mathcal{K}^a_{21}$ operator, which
    is a compact perturbation of $\frac{\beta}{2}\frac{\dd }{\dd s}\mathcal{H}$,
    and the jump in the $\mathcal{K}_{21}^b \mathcal{H}$ operator.
    The jump properties of $\mathcal{K}_{21}^b \mathcal{H}$ depend on the side of approach, while the singularity of
    $\mathcal{K}_{21}^a$ does not. Relatedly, the integral in the 
    on-surface operator for $\mathcal{K}^a_{21}$ must be understood 
    in the Hadamard finite part sense. Fortunately, owing to the cancellation 
    that occurs, this integral never needs to be evaluated directly.
\end{rmk}
}
{ 
\begin{thm}
\label{thm:free}
    Let $\Omega$ be a simply-connected domain. 
    Let $\cK_1=\cK_1^a + \beta^{\pm} \cK_1^b\mathcal{H}$ and
    let $\cK_1^a$, $\cK_1^b$, and $\cK_2$ be the layer potentials corresponding to the
    kernels $K^a_1$, $K_1^b$, and $K_2$, respectively. Let 
    $\cK_{11} = \cK_{11}^a + \beta^\pm \cK_{11}^b \mathcal{H}$ 
    and $\cK_{21} = \cK_{21}^a + \beta^\pm \cK_{21}^b\mathcal{H}$. Let 
    $\cK_{11}^a,\cK_{11}^b,\cK_{12},\cK_{21}^a,\cK_{21}^b,$ and 
    $\cK_{22}$ 
    be the boundary operators corresponding to the kernels $K_{11}^a,K_{11}^b,K_{12},K^a_{21},K_{21}^b$, 
    and $K_{22}$, respectively. For ease of exposition, we set
    \begin{align} 
    \mathcal{B}_1(\bx_0) &= \nu \Delta + (1-\nu) \left (\bn(\bx_0) \cdot \nabla_{\bx}\right)^2 \, ,\\
    \mathcal{B}_2(\bx_0) &= \left (\bn(\bx_0) \cdot \nabla_{\bx}\right)^3 + (2-\nu) \left (\btau(\bx_0) \cdot \nabla_{\bx}\right)^2 (\bn(\bx_0)\cdot \nabla_\bx) \nonumber \\ &\qquad \qquad + (1-\nu)\kappa(\bx_0) \left ( 
    (\btau(\bx_0) \cdot \nabla_\bx)^2 - (\bn(\bx_0)\cdot \nabla_\bx)^2 \right )  \, ,
    \end{align}
    which are the differential operators appearing in the free plate boundary 
    conditions,~\eqref{freebcs}. After inserting the integral operators defined above into the boundary conditions, we obtain the following jump relations:
\begin{align}
    \lim_{\bx \to \bx_0^{\pm}} \mathcal{B}_1(\bx_0) \mathcal{K}_{1}[\sigma](\bx) &= \mp \frac{1}{2} \sigma(\bx_0) + \mathcal{K}_{11}[\sigma](\bx_0) \, , \\
    \lim_{\bx \to \bx_0^{\pm}} \mathcal{B}_1(\bx_0) \mathcal{K}_{2}[\sigma](\bx) &= \mathcal{K}_{12}[\sigma](\bx_0) \, , \\
    \lim_{\bx \to \bx_0^{\pm}} \mathcal{B}_2(\bx_0) \mathcal{K}_{1}[\sigma](\bx) &= -\frac{\beta}{2} \frac{\dd}{\dd s} \mathcal{H}[\sigma](\bx_0) + \mathcal{K}_{21}[\sigma](\bx_0) \, , \\
    \lim_{\bx \to \bx_0^{\pm}} \mathcal{B}_2(\bx_0) \mathcal{K}_{2}[\sigma](\bx) &= \pm \frac{1}{2} \sigma(\bx_0) + \mathcal{K}_{22}[\sigma](\bx_0) \, ,
\end{align}
for $\beta = \displaystyle\frac{1+\nu}{2}$ and $\bx_0 \in \partial\Omega$, where $\bx_0^+$ corresponds to the limit from the exterior and $\bx_0^-$ corresponds to the limit from the interior, each along a normal direction.

Moreover, the corresponding exterior and interior integral equations
for the integral representation $u=\cK_1[\rho_1] + \cK_2[\rho_2]$, i.e. 
\begin{equation}
    \begin{pmatrix}
        \mp \frac{1}{2} I & 0\\
        -\frac{\beta}{2} \frac{\dd}{\dd s} \mathcal{H} & \pm\frac{1}{2} I
    \end{pmatrix} \begin{pmatrix} \rho_1 \\ \rho_2 \end{pmatrix} + \begin{pmatrix} 
    \cK_{11} & \cK_{12} \\ 
    \cK_{21} & \cK_{22}
    \end{pmatrix} \begin{pmatrix} \rho_1 \\ \rho_2 \end{pmatrix} = \begin{pmatrix} f_1 \\ f_2 \end{pmatrix} \, ,
\end{equation}
are second kind on $L^2(\partial \Omega) \times L^2(\partial \Omega)$ for
$\beta^2 \ne  1$. These equations have the equivalent form 
\begin{equation} \label{eq:freebie2}
    \begin{pmatrix}
        \left (\mp \frac{1}{2} \pm \frac{\beta^2}{2} \right) I & 0\\
         0 & \pm\frac{1}{2} I
    \end{pmatrix} \begin{pmatrix} \rho_1 \\ \rho_2 \end{pmatrix} + \begin{pmatrix} 
    \cK_{11}^a + \beta^\pm\left (\cK_{11}^b +\frac{\beta}{2} \mathcal{H}\right )\mathcal{H} - 2\beta\beta^\pm \mathcal{D}^2 & \cK_{12} \\ 
    \left ( -\frac{\beta}{2} \frac{\dd}{\dd s} \mathcal{H} + \cK^a_{21}\right )  + \cK_{21}^b \mathcal{H} & \cK_{22}
    \end{pmatrix} \begin{pmatrix} \rho_1 \\ \rho_2 \end{pmatrix} = \begin{pmatrix} f_1 \\ f_2 \end{pmatrix} \, .
\end{equation} 
\end{thm}}

{A similar integral representation results in a SKIE for the multiply-connected case.
We omit the formula for the integral equation system in this case for 
the sake of brevity.

\begin{cor}
    Let $\Omega = \bigcup_{i=1}^M\Omega_i$ where the $\Omega_i$ are 
    bounded, simply connected domains with smooth boundaries and disjoint
    closures. Let 
    $\cK_1^{a,\partial \Omega_i}$, $\cK_{1}^{b,\partial \Omega_i}$,
    and $\cK_{2}^{\partial \Omega_i}$ denote the layer potential operators 
    for the kernels $K_1^a$, $K_1^b$, and $K_2$, respectively, with 
    sources restricted to $\partial \Omega_i$. Let $\mathcal{H}^{\partial \Omega_i}$ 
    denote the Hilbert transform operator on $\partial \Omega_i$. 
    Finally, let $\rho_1^{\partial \Omega_i}$ and $\rho_2^{\partial \Omega_i}$ 
    denote the restrictions of densities $\rho_1$ and $\rho_2$ to 
    $\partial \Omega_i$. If we represent
    the solution $u$ of the free plate problem using the layer potential 
    representation 

\begin{equation}
     u = \sum_{i=1}^M \left (\cK_1^{a,\partial \Omega_i} \left [\rho_1^{\partial \Omega_i} \right] 
     + \cK_1^{b,\partial \Omega_i} \left [\mathcal{H}^{\partial \Omega_i}\left [\rho_1^{\partial \Omega_i}\right ] \right] 
     + \cK_2^{\partial \Omega_i} \left [\rho_2^{\partial \Omega_i} \right] \right ) \, ,
\end{equation}    
then the resulting system of boundary integral equations is second kind.
\end{cor}

Before we prove Theorem~\ref{thm:free}, we require some preliminary results.}
{We will use the 
following relation between $\mathcal{H}$ and 
$\mathcal{D}$}, which is a special case of the Poincar\'e-Bertrand formula; see, 
for example, \cite{muskhelishvilisingularbook,shidong} and the references therein. 
\begin{lem} \label{lem:pb}
{Assume that $\Omega$ is a bounded, simply connected domain. 
Let $\mathcal{H}$ and $\mathcal{D}$ be as defined above. Then,}
\begin{equation}
    \frac{1}{4}\mathcal{H}^2 = - \frac{I}{4} +  \mathcal{D}^2 \, . \label{hilbsquared} 
\end{equation}
\end{lem}

Let $K^{\mathcal{H}}$ be the kernel of the Hilbert transform \eqref{hilbertdef} and let $K^{\mathcal{H}'} := \frac{\dd }{\dd s} K^{\mathcal{H}}$ be the kernel of {$\frac{\dd}{\dd s} \mathcal{H}$} (given below).
We have the following result:

\begin{lem}
{The kernels $K^a_{11},(K^b_{11} + \frac{\beta}{2} K^{\mathcal{H}} ),(K^a_{21} - \frac{\beta}{2} K^{\mathcal{H}'}), K^b_{21}$, and $K_{22}$ are continuous, while the kernel $K_{12}$ is weakly (log) singular when restricted to the boundary}. \label{freesmooth}
\end{lem} 

\begin{proof}
As before, {it is} sufficient for us to show the continuity of the corresponding kernels {with $G$ replaced by} the biharmonic Green's function $G^{B}(\bx,\by)$. Suppose $\bgamma(t)$ is the arc length parametrization {of} the boundary $\partial \Omega$. Using the Taylor expansions {\eqref{gammaseries}-\eqref{nseries}} together with the formulae in \ref{appfree}, we have the following asymptotics for the terms in $K_{11}^a$:
    \begin{align}
        G^{B}_{\bn_\bx  \bn_\bx  \bn_\by }{(\bgamma(0),\bgamma(s))} &= - \frac{1}{8\pi}\kappa(0) + \mathcal{O}(s) \, , \label{gnxnxnyasymp} \\
        G^{B}_{\btau_\bx  \btau_\bx  \bn_\by } {(\bgamma(0),\bgamma(s))}  &= \frac{3}{8\pi}\kappa(0) + \mathcal{O}(s) \,  . \label{gtxtxnyasymp}
    \end{align}
{From this we obtain the following on-surface limit:}
    \begin{equation}
        \lim_{s\to 0}  K^a_{11}{(\bgamma(t) ,\bgamma(t+s) )} = \frac{3\nu-1}{8 \pi}\kappa(t) \, .
    \end{equation}
    For a {sufficiently} smooth curve, this limit is well-defined, so we conclude that $K_{11}^a$ is continuous. We now analyze $(K^b_{11} + \frac{\beta}{2} K^{\mathcal{H}} )$. Substituting the expansions \eqref{gammaseries}-\eqref{nseries} into the formulae for $G^{B}_{\bn_\bx  \bn_\bx  \btau_\by } , G^{B}_{\btau_\bx  \btau_\bx  \btau_\by }$, and $K^{\mathcal{H}}$, we obtain: 
    \begin{align}
         G^{B}_{\bn_\bx  \bn_\bx  \btau_\by } {(\bgamma(0),\bgamma(s))}  &= \frac{1}{4\pi} \frac{1}{s} + \mathcal{O}(s) \, , \label{gnxnxtyasymp} \\
          G^{B}_{\btau_\bx  \btau_\bx  \btau_\by } {(\bgamma(0),\bgamma(s))} &= \frac{1 }{4\pi} \frac{1}{s} + \mathcal{O}(s) \, , \label{gtxtxtyasymp} \\
         K^{\mathcal{H}} {(\bgamma(0),\bgamma(s))} &= -\frac{1}{\pi} \frac{1}{s} + \mathcal{O}(s) \, , \label{hilbexpansion}
    \end{align}
    which gives:
    \begin{equation}
        \lim_{s\to 0} \left[ K^b_{11}{(\bgamma(t) ,\bgamma(t+s) )} + \frac{\beta}{2} K^{\mathcal{H}}{(\bgamma(t) ,\bgamma(t+s) )} \right] = 0 \, ,
    \end{equation}
    {thereby making the kernel $(K^b_{11} + \frac{\beta}{2} K^{\mathcal{H}} )$ continuous. We now turn to $(K^a_{21}- \frac{\beta}{2} K^{\mathcal{H}'})$. The formula for $ K^{\mathcal{H}'}$ is obtained by taking one tangential derivative of $K^\mathcal{H}$ with respect to the target variable:}
    \begin{equation}
        K^{\mathcal{H}'}(\bx ,\by ) = \frac{1}{\pi} \frac{ \btau (\by)  \cdot \btau (\bx)}{\|\bx -\by \|^2} - \frac{2}{\pi} \frac{[(\bx -\by )\cdot \btau (\bx)][(\bx -\by ) \cdot  \btau (\by) ]}{\|\bx -\by \|^4} \, .
    \end{equation}
    As before, {we} substitute the Taylor expansions into the formulae {for} $G^{B}_{\bn_\bx  \bn_\bx  \bn_\bx  \bn_\by }$, $G^{B}_{\bn_\bx  \btau_\bx \btau_\bx  \bn_\by }$, and $K^{\mathcal{H}'}$, which yields: 
    \begin{align}
        G^{B}_{\bn_\bx  \bn_\bx  \bn_\bx  \bn_\by }{(\bgamma(0),\bgamma(s))} &= -\frac{3}{4\pi}\frac{1}{s^2} + \frac{5}{16\pi}\kappa^2(0) + \mathcal{O}(s) \, , \\
        G^{B}_{\bn_\bx  \btau_\bx \btau_\bx  \bn_\by } {(\bgamma(0),\bgamma(s))} &= \frac{1}{4\pi} \frac{1}{s^2} - \frac{17}{48 \pi} \kappa^2(0) + \mathcal{O}(s) \, , \\
        K^{\mathcal{H'}}{(\bgamma(0),\bgamma(s))} &= -\frac{1}{\pi} \frac{1}{s^2}  - \frac{1}{12\pi}\kappa^2(0)+ \mathcal{O}(s) \, .
    \end{align}
    Combining the formulae above with \eqref{gnxnxnyasymp}-\eqref{gtxtxnyasymp} yields the limit:
    \begin{equation}
        \lim_{s\to 0} \left[ K^a_{21} (\bgamma(t) ,\bgamma(t+s) ) - \frac{\beta}{2} K^{\mathcal{H}'} (\bgamma(t) ,\bgamma(t+s) ) \right]  = \frac{1-\nu}{8\pi} \kappa^2(t ) \, .
    \end{equation}
    Since $\partial \Omega $ is {sufficiently} smooth, we conclude that the kernel {$(K^a_{21}- \frac{\beta}{2} K^{\mathcal{H}'})$} is indeed continuous. 
    
    Now we examine $K_{21}^b$. Substituting the Taylor expansions into the closed-form formulae for $G^{B}_{\bn_\bx  \bn_\bx  \bn_\bx \btau_\by  }$ and $G^{B}_{\btau_\bx  \btau_\bx  \btau_\by  \bn_\bx }$, we get: 
    \begin{align}
        G^{B}_{\bn_\bx  \bn_\bx  \bn_\bx \btau_\by  } {(\bgamma(0),\bgamma(s))} &= \frac{1}{8 \pi} \kappa'(0) + \mathcal{O}(s) \, , \\
        G^{B}_{\bn_\bx \btau_\bx  \btau_\bx  \btau_\by  }{(\bgamma(0),\bgamma(s))} &= -\frac{1}{24\pi}\kappa'(0) +  \mathcal{O}(s) \, .
    \end{align}
    Substituting these expressions as well as \eqref{gnxnxtyasymp}-\eqref{gtxtxtyasymp} we get:
    \begin{equation}
        \lim_{s\to 0} K^b_{21} {(\bgamma(t) ,\bgamma(t+s) )} = \beta \frac{1+\nu}{24\pi}\kappa'(t ) \, .
    \end{equation}
     Since $\partial \Omega $ is {sufficiently} smooth, we conclude that $K^b_{21}$ is also continuous.
     
     Finally, we analyze $K_{22}$. Taylor expanding the terms $G^B_{\bn_\bx \bn_\bx \bn_\bx}$ and $G^B_{\bn_\bx \btau_\bx \btau_\bx}$ we have: 
     \begin{align}
         G^B_{\bn_\bx \bn_\bx \bn_\bx}{(\bgamma(0),\bgamma(s))} &= \frac{3}{8\pi} \kappa(0) + \mathcal{O}(s) \, , \\
         G^B_{\bn_\bx \btau_\bx \btau_\bx}{(\bgamma(0),\bgamma(s))} &= - \frac{1}{8\pi}\kappa(0) + \mathcal{O}(s) \, ,
     \end{align}
     Meanwhile, the term $ G^{B}_{\btau_\bx  \btau_\bx } -  G^{B}_{\bn_\bx  \bn_\bx }$ {simplifies to}: 
     \begin{equation}
         G^{B}_{\btau_\bx  \btau_\bx }(\bx,\by) -  G^{B}_{\bn_\bx  \bn_\bx }(\bx,\by) = \frac{1}{4\pi} - \frac{1}{2\pi} \frac{[(\bx -\by ) \cdot  \bn (\bx ) ]^2}{\|\bx -\by \|^2} \, .
     \end{equation}
     {Again considering Taylor approximations of the kernels near $s=0,$ it follows immediately that the second term is $\mathcal{O}(s^2)$.} Therefore: 
     \begin{equation}
         \lim_{s\to  0} K_{22} {(\bgamma(t) ,\bgamma(t+s) )} = \frac{3-\nu}{8} \kappa(t) \, ,
     \end{equation}
    which is again continuous provided sufficient regularity of the boundary. Meanwhile, the terms in $K_{12}$ have {only}
    two derivatives so that they are {logarithmically} singular, which can be seen directly {from} the formulae in \ref{appfree}. 
\end{proof} 

\begin{proof}[Proof of Theorem~\ref{thm:free}.] 
{The jump properties are established in \ref{app:jumpcondfree}.
The equivalent formulation, \eqref{eq:freebie2}, follows by applying 
Lemma~\ref{lem:pb} to the formula 
\begin{equation}
    \mathcal{K}_{11}^b \mathcal{H} = \left ( \cK_{11}^b + \frac{\beta}{2} \mathcal{H}\right ) 
\mathcal{H} - \frac{\beta}{2} \mathcal{H}\mathcal{H} \, .
\end{equation} }

{It remains to show that the integral operators $\mathcal{D}^2$, $\cK_{11}^a$,
$(\mathcal{K}_{11} + \tfrac{\beta}{2}\mathcal{H})\mathcal{H}$ , 
$\mathcal{K}_{12}$, $-\frac{\beta}{2} \frac{\dd}{\dd s} \mathcal{H} + \mathcal{K}_{21}$, 
and $\mathcal{K}_{22}$ are compact.  }
{It is known from the theory of integral equations \cite{kress} 
that integral operators on $\partial \Omega$ with continuous or weakly singular kernels are compact from $L^2(\partial \Omega) \to L^2(\partial \Omega)$.} In particular, since the kernels are continuous for 
$K^a_{11},(K^b_{11} + \frac{\beta^2}{2} K^{\mathcal{H}} ),(K^a_{21}- \frac{\beta}{2} K^{\mathcal{H}'}), K^b_{21}$, and $K_{22}$, and logarithmically singular in the case of $K_{12}$, they correspond to compact integral operators. Furthermore, the Hilbert transform is a {bounded operator on $L^2(\partial \Omega)$ (see, {\em inter alia}, \cite{duoandikoetxea2001fourier}),} and since the composition of a bounded operator with a compact operator is again compact, we conclude that $(\mathcal{K}_{11}^b+\tfrac{\beta}{2} \mathcal{H}) \mathcal{H}$ and $\mathcal{K}_{21}^b \mathcal{H}$ are compact. {Finally, it is well known that 
$\mathcal{D}$ is compact on a smooth boundary.}
\end{proof}

{
\begin{rmk}
    In practice, we discretize the equivalent formulation,
    \eqref{eq:freebie2}, because the cancellations in the operators
    are more explicit.
\end{rmk} }

\section{Numerical implementation}\label{num_imp}

In this section, we sketch relevant details for the numerical implementation of the proposed scheme
for solving the boundary value problems. In subsection~\ref{num_details}, we briefly discuss the discretization approach we use. Following this, in subsection~\ref{num_cancel} we touch on catastrophic cancellations which arise in a na\"{i}ve attempt to discretize the integral operators, as well as a method for circumventing these difficulties.

\subsection{Details of the discretization}
\label{num_details}
In order to solve boundary integral equations  of the form \eqref{eq:biesystem}, we use a standard modified Nystr\"{o}m method. In particular,
the boundary curve $\partial \Omega$ is divided into $N_p$ curved panels, each of which
is approximated by a polynomial interpolant at scaled, $n_{\textrm{GL}}$-th order
Gauss-Legendre nodes, for a total of $N=N_p n_{\textrm{GL}}$ boundary
nodes; see Figure~\ref{convergence_figure} for an illustration of the 
discretization nodes with $N_p=16$ and $n_{\textrm{GL}}=16$. Unless otherwise stated,
the numerical results are for $n_{\textrm{GL}}=16$. We approximate the boundary layer
densities $\rho_1,\rho_2$ by their polynomial interpolants at the same nodes. 

The discrete approximation of \eqref{eq:biesystem} is then obtained by enforcing
that the integral equation holds at the boundary nodes. 
Since the integral kernels $K_{ij}$ are in general continuous but not smooth, the integrals in 
\eqref{eq:biesystem} require special quadrature rules to obtain a high order accurate
discretization. 
Here, we apply special generalized Gaussian quadrature rules~\cite{bremer2010nonlinear} (GGQ) which are 
pre-computed for the class of integral kernels we encounter, specifically integral 
kernels of the form $\psi(t) + \phi(t) \log |t-t_0|$, where $\psi$ and $\phi$ 
are smooth functions and $t_0$ is a given point. 

To visualize the solutions we must also evaluate the integrals in the layer 
potential representation, \eqref{eq:layerpotential}. For target points, $\bx $,
near the boundary, the integral is nearly-singular. In this case, we approximate the corresponding 
integrals with high accuracy using adaptive Gaussian integration. For the numerical results presented below, we use the \texttt{chunkIE} package in MATLAB
\cite{Askham_chunkIE_a_MATLAB_2024}, which provides utilities for discretizing the curve 
by panels, applying generalized Gaussian quadrature rules to obtain the 
discretization of \eqref{eq:biesystem}, and computing the integrals
in the layer potential representation by adaptive quadrature. 

\subsection{Catastrophic cancellation}
\label{num_cancel}

While the integral kernels derived in this work are at most weakly
singular when restricted to the boundary, a na\"{i}ve implementation
of the integral kernel formulae can lead to catastrophic cancellation. Recall that the Green's function is defined as a scaled difference of the Helmholtz
and modified Helmholtz Green's functions, i.e.
$$ G(\bx ,\by ) = \frac{1}{2k^2} \left ( \frac{i}{4} H_0^{(1)}(k\|\bx -
\by \|) - \frac{1}{2\pi} K_0(k\|\bx - \by\|) \right ) \, .$$
However, $G$ and its derivatives should not be evaluated by evaluating the
$H_0^{(1)}$ and $K_0$ terms separately and taking the difference for
small values of $|k|\|\bx - \by\|$. These terms have the same leading
order asymptotics for small $|k|\|\bx - \by\|$, so $G$ and its derivatives
should be evaluated using the appropriate Taylor series~\cite[\S 10.8]{NIST}
and explicitly performing the cancellation between these terms, as in
\eqref{eq:gseries}. In fact, it is convenient to use these power series
to write a subroutine which can also evaluate $G(\bx,\by) - G^B(\bx,\by)$
stably for small $|k|\|\bx-\by\|$, as we will see below. 

Other forms of numerical cancellation can occur in implementing the
integral kernels, which are given as linear combinations of the derivatives
of $G$. Indeed, such cancellations are necessary to avoid singularities
in the integral kernels. Conveniently, up to $5$th order derivatives of the
quantity $G(\bx,\by)-G^B(\bx,\by)$ are bounded and continuous. {The dominant
numerical cancellations that arise from taking linear combinations of derivatives
of the Green's function} can then be mitigated by carefully considering the
same linear combinations of $G^B$. {This process is more straightforward
for kernels that only involve directional derivatives of the Green's function
and can be more complicated for kernels that involve the curvature and its
derivatives.}

For example, consider the kernel $K_{21}$ from the clamped plate boundary
integral equation:

$$  K_{21} = G_{\bn_\bx  \bn_\by  \bn_\by  \bn_\by } +
3G_{\bn_\bx  \bn_\by  \btau_\by  \btau_\by } \, .$$
This kernel is shown to be continuous above but the individual terms are
proportional to $1/\|\bx-\by\|^2$. Fortunately, in the biharmonic part $K_{21}^B$, the cancellations can be carried out explicitly. We have
$$ G^B_{\bn_\bx  \bn_\by  \bn_\by  \bn_\by } +
3G^B_{\bn_\bx  \bn_\by  \btau_\by  \btau_\by } =
\frac{4}{\pi} \frac{[\br \cdot \bn (\by)]^3 [\br \cdot \bn (\bx)]}{\|\br\|^6}
- \frac{3}{\pi} \frac{[\br \cdot \bn (\by)]^2 [\bn(\by)\cdot \bn(\bx)]}{\|\br\|^4} \, ,$$
where $\br = \bx-\by$. There is a milder
form of numerical cancellation in evaluating $\br\cdot \bn (\by) =
\mathcal{O}(\|\br\|^2)$, but this did not appear to have a significant effect
in the numerical results. A reasonably stable scheme for evaluating
$K_{21}$ when $|k|\|\bx-\by\|$ is small is then obtained by
using {a series expansion} to evaluate $(G-G^B)_{\bn_\bx  \bn_\by  \bn_\by  \bn_\by }
+ 3(G-G^B)_{\bn_\bx  \bn_\by  \btau_\by  \btau_\by }$, using a
straightforward implementation to evaluate
$G^B_{\bn_\bx  \bn_\by  \bn_\by  \bn_\by } +
3G^B_{\bn_\bx  \bn_\by  \btau_\by  \btau_\by }$
via the formula above, and then adding the results.

In some cases it is more difficult to carry out the cancellations in the biharmonic part. Consider the kernel $K_{21}$ for the supported plate boundary integral equation:
\begin{multline*}
K_{21} = G_{\bn_\bx  \bn_\bx  \bn_\by  \bn_\by  \bn_\by } + \nu  G_{\btau_\bx  \btau_\bx  \bn_\by  \bn_\by  \bn_\by }   + \alpha_1 G_{ \bn_\bx  \bn_\bx  \bn_\by  \btau_\by  \btau_\by } + \nu  \alpha_1 G_{\btau_\bx  \btau_\bx  \bn_\by  \btau_\by  \btau_\by } +  \\
   \alpha_2 \kappa(\by ) G_{\bn_\bx  \bn_\bx  \bn_\by  \bn_\by  } +\nu \alpha_2 \kappa(\by ) G_{\btau_\bx  \btau_\bx  \bn_\by  \bn_\by  } + \alpha_3 \kappa'(\by ) G_{\bn_\bx  \bn_\bx  \btau_\by  } + \nu \alpha_3 \kappa'(\by ) G_{\btau_\bx  \btau_\bx  \btau_\by  }\, .
\end{multline*}
A proper treatment of the numerical cancellations in this kernel for small
$\|\bx-\by\|$ would require constructing a higher order Taylor approximation of the coordinates 
of the curve. In lieu of this, we instead leverage the smoothness of the biharmonic
part of $K_{21}$ to (partially) mitigate such cancellations. The idea is to split the kernel into $K_{21} = (K_{21}-K_{21}^B) + K_{21}^B$
as before. For the difference, $(K_{21}-K_{21}^B)$, the power series
for $G-G^B$ can be used to stably evaluate the kernel. This component
of the kernel is continuous, with limited smoothness owing to the logarithmic
terms in the power series. We thus treat this part of the term using
a GGQ rule for logarithmic singularities.

To handle $K_{21}^B$, we observe that the use of a GGQ rule is unnecessary
because the kernel is smooth (on a smooth curve). We can instead use a
standard Gaussian quadrature rule. This provides a significant improvement
in stability which can be understood by analyzing the inter-node spacing
and integral weights in each scheme.
Consider Gauss-Legendre nodes of order $n_{\textrm{GL}}=16$
scaled to $[-1,1]$. The smallest inter-node spacing occurs for the
node nearest $-1$, at a distance $\approx 4.5 \times 10^{-2}$ and
the corresponding (smooth) integration weight is $\approx 6.2 \times 10^{-2}$.
In contrast, the GGQ rule we apply for functions with logarithmic singularities
uses support nodes which are specific to each ``target'' Gauss-Legendre node.
For the node of order $n_{\textrm{GL}}=16$ that is closest to $-1$, the
closest support node is at a distance of $\approx 2.2 \times 10^{-4}$
with a weight of $\approx 2.4 \times 10^{-3}$. Because we expect to lose
precision proportional to the product of the integral weight and
the inverse square of the inter-node distance (when using an unstable
evaluator for $K_{21}^B$), we expect that the standard rule will yield 
higher precision and this is observed in practice. 

\section{Results} \label{results}

In order to test the error of these methods, we place a point source {inside $\Omega$} and set the righthand side of the integral equation to be the {boundary data corresponding to the point source}. Because the field resulting from this point source satisfies the {time harmonic} flexural wave equation {outside of $\Omega$}, the integral equation should reconstruct the Green's function {in the exterior, $E$}. We refer to this test as the analytic solution test, since the analytic expression for the Green's function is known. 

\begin{figure}[ht]
\centering
\includegraphics[scale=1.3]{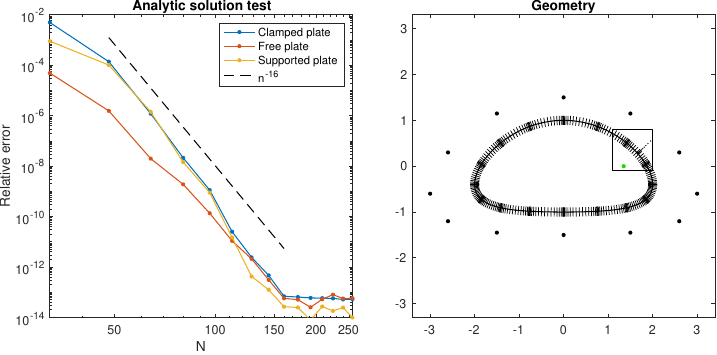}
\caption{Convergence of the three boundary integral equations as a function of the number of discretization points {$N$ ($k = 8, \nu = 1/3$)}. The solutions converge at sixteenth order (left). The {equations} were solved on {the boundary of a} droplet (right), with the green dot representing the location of the point source and the black dots representing the locations where the error was measured. { The discretization plotted on the right was the finest one tested ($N_p=16, n_{\textrm{GL}}=16$). The black box represents the region near the boundary that was chosen for closer investigation (see Figure \ref{near_boundary_figure}).}}\label{convergence_figure}
\end{figure}

\begin{figure}[ht]
\centering
\includegraphics[scale=1.3]{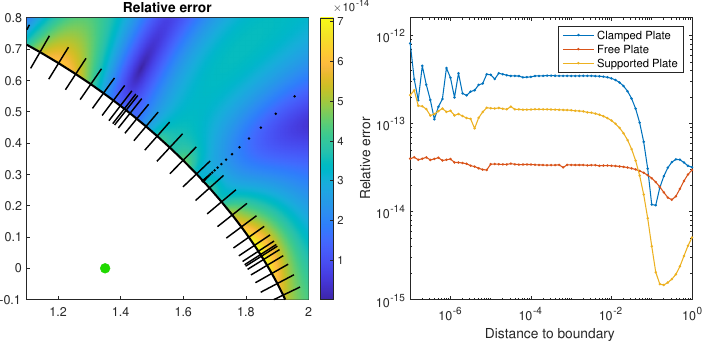}
\caption{{The accuracy of evaluation for points near the boundary was checked for the same geometry and discretization as Figure \ref{convergence_figure}. The relative error for the free plate problem is plotted in some region near the boundary (left). The error was also calculated at a set of points that get exponentially close to the boundary, indicated by black dots, and plotted for the clamped, free, and supported plate problems (right). }}\label{near_boundary_figure}
\end{figure}

{In particular,} we solve the BIE {on the boundary of a droplet} (Figure \ref{convergence_figure}), which is parameterized by the equations:
\begin{align}
x(t) &= 2 \cos(t) \, , \\
y(t) &= \sin(t) - 0.4 \cos^2(t) \, ,
\end{align}
with {$t \in [-\pi, \pi)$.} The source was placed {inside $\Omega$} at the point $(x,y) = (1.35, 0).$ The error was measured at a collection of twelve points that come from the {rescaled curve} $\frac{3}{2} ( x(t), y(t) )$, where $t = n \pi / 6$ for {$n = -6, ..., 5$.} To compute the relative error, the $\ell_\infty$ error at the collection of points is calculated and then divided by the sums of the $L_1$ norms of the densities on the boundary. 
{We also plot the error in the computed solution using the finest grid ($N=256$)
for clamped, supported, and free boundary conditions as the evaluation point approaches the boundary in Figure~\ref{near_boundary_figure}.
This is a more stringent test of the precision in the computed density.}

{

\begin{figure}[ht]
\centering
\includegraphics[scale=1.2]{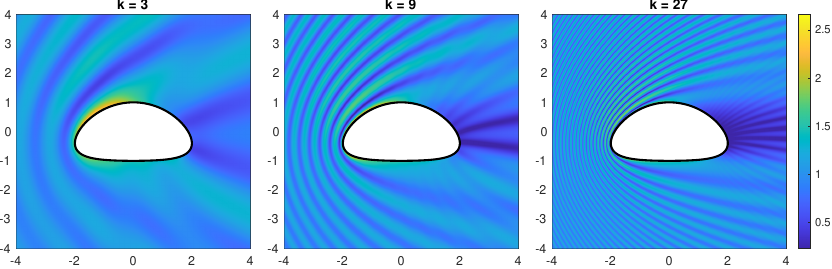}
\caption{{Wave scattering off a droplet with free plate boundary conditions for three different wavenumbers. The field is incident from the left and the magnitude of the total field $|u|$ is plotted. }}\label{variable_k_figure}
\end{figure} 

The next two figures demonstrate how the qualitative behavior of wave scattering for the free plate problem depends on the wavenumber $k$ (Figure~\ref{variable_k_figure}) and Poisson's ratio $\nu$ (Figure~\ref{variable_nu_figure}). To solve the wave scattering problem, we write the total field as a sum of incident and scattered fields $u(\bx) = u^{(i)}(\bx) + u^{(s)}(\bx)$, where the incident field $u^{(i)}(\bx) = e^{i \mathbf{k} \cdot \bx}$ is a plane wave. Then, we let the righthand side of the BIE be the negative of the boundary data corresponding to a plane wave with $\|\mathbf{k}\| = k$. 

\begin{figure}[ht]
\centering
\includegraphics[scale=1.2]{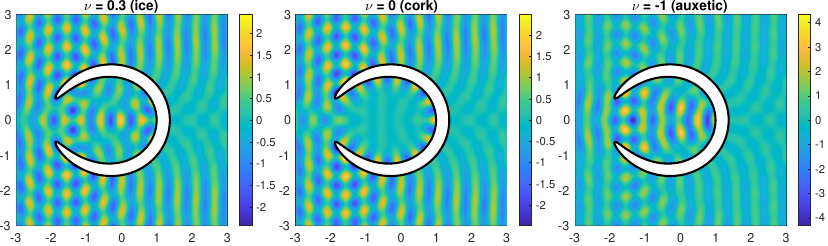}
\caption{{Wave scattering in free plates with three different Poisson's ratios: $\nu = 0.3$ (ice), $\nu = 0$ (cork), and $\nu = -1$ (auxetic). The field is incident from the left and the wavenumber is $k = 12$. The quantity that is plotted is the real part of the total field. The condition numbers for these three cases were 1.37E+04 ($\nu = 0.3$), 1.90E+04 ($\nu = 0$), and 2.46E+04 ($\nu = -1$). For concave domains, Poisson's ratio has a large effect on wave reflection and resonant behavior.
}}\label{variable_nu_figure}
\end{figure} 

In the case of variable wavenumber (Figure~\ref{variable_k_figure}), low wavenumbers lead to larger amplitudes of the solution near the boundary and a slightly muddied appearance of the total field, while larger wavenumbers lead to a darker shadow and a more distinct wake. In order to test the effect of Poisson's ratio $\nu$ in the boundary conditions \eqref{freebcs}, the free plate problem was solved for three different materials: ice ($\nu = 0.3$), cork ($\nu = 0$), and auxetic ($\nu = -1$).  Note that a concave geometry was chosen to exaggerate the differences in the total field. Different Poisson's ratios lead to different resonant patterns and amplitudes inside this shape (Figure~\ref{variable_nu_figure}). When the geometry is convex, there is a much more subtle dependence of the total field on Poisson's ratio. }

\begin{figure}[ht]
\includegraphics[scale=1.24]{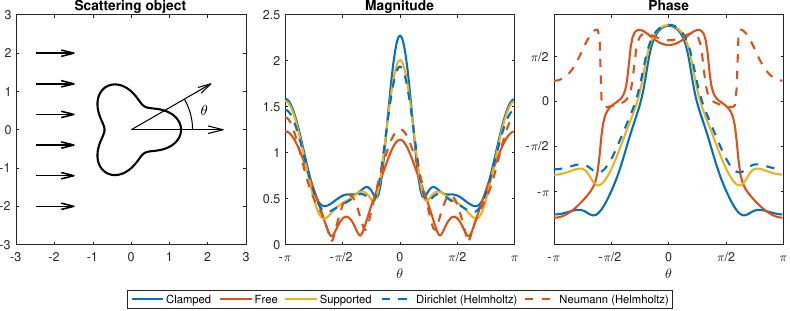}
\caption{Comparison of {the scattered field} for various boundary conditions {($k = 3, \nu = 1/3$), measured in the far field}. {The geometry is a starfish (left) with arrows indicating the direction of the incident field and the angle $\theta$}. Different boundary conditions lead to differences {in both magnitude (middle) and phase (right) of the scattered field }. }\label{farfield1}
\end{figure}

To investigate the qualitative differences {between} boundary conditions, we analyze the far field patterns that are generated when a wave is scattered by an object, in this case, a tripedal starfish (Figure \ref{farfield1}). The {general form for the parametrization of a `starfish' domain} is given by
\begin{align}
    x(t) &= (1+ A \cos(n t)) \cos(t) \, , \\
    y(t) &= (1+ A \cos(n t)) \sin(t) \, ,
\end{align}
where {$t \in [-\pi, \pi)$,} $A$ is the size of the arms, and $n$ is the number of appendages. Because the Yukawa part of the Green's function, {i.e. the $K_0$ term in \eqref{greens},} decays exponentially away from the boundary, the solutions far away from the boundary will naturally look like solutions to the Helmholtz equation. The limiting behavior of the {Hankel part of the Green's function is given by} $H_0^{(1)}(k\|\bx \|) \sim \displaystyle \sqrt{\frac{2}{\pi k \|\bx  \| }} e^{i ( k \|\bx  \| - \frac{1}{4} \pi )} $ \cite{NIST}. { Thus, it is possible to show that the scattered field $u^{(s)}(\bx)$ has similar $ e^{i  k \|\bx  \|} /\sqrt{\|\bx\|}$ decay to leading order. As a measure of the far field, we send a plane wave with wavenumber $k = 3$ toward the scattering object and solve for the scattered field $u^{(s)}(\bx)$ on a large ball of radius $\| \bx \| = 1000$ (roughly 480 wavelengths).  We plot both the magnitude and phase of $f(\theta) :=  \sqrt{\|\bx  \|} e^{- i k \| \bx  \|} u^{(s)}(\bx)$, where the phase is given by the formula {$\phi(\theta) := \arctan \left( \Im(f(\theta)) / \Re(f(\theta)) \right)$}.} To understand how the solutions of flexural wave problems might differ from the solutions to the Helmholtz equation, we also plot {the far field patterns for the Helmholtz scattering problem} with Dirichlet and Neumann boundary conditions, { for comparison. }

\begin{figure}[ht]
\includegraphics[scale=1.24]{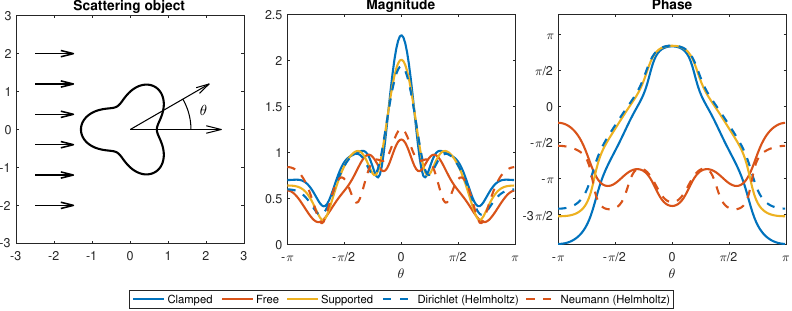}
\caption{Far field patterns generated by the same starfish {as Figure \ref{farfield1},} rotated by $\pi/3${, for the same wavenumber and Poisson's ratio}. The change in the orientation in the starfish {(left) leads to slightly different magnitude (middle) and phase (right) of the scattered field}. }\label{farfield2}
\end{figure}

{Unsurprisingly, the choice of} boundary condition leads to very different patterns in the far field (Figures \ref{farfield1} and \ref{farfield2}). The clamped plate has a disruptive effect on the incoming wave, leading to the darkest shadow and the strongest backscatter. This effect is stronger than the analogous Dirichlet problem for the Helmholtz equation. Meanwhile, the free plate has the lightest shadow and weakest backscatter, suggesting that the wave is able to pass through the object with less disturbance. Again, this passivity is stronger in the flexural wave case than in the analogous Neumann problem for the Helmholtz equation. The far field of the supported plate is similar to that of the Dirichlet problem for the Helmholtz equation. { These findings agree with the qualitative descriptions of wave dynamics for these problems found in classic texts on thin plate theory \cite{landau59,timoshenko1959theory} }.

The phase of the far field also varies to great extent. In some parts of the far field ($\theta = \pm {\pi/4}$), the clamped and free plate are roughly half a period out of phase, while in other parts of the far field ($\theta = \pi$), the clamped and far fields are {roughly} in phase. This is in contrast with the Helmholtz case, where the Dirichlet and Neumann problems are once again half a period out of phase at the angle $\theta = \pi$.

\begin{figure}[ht]
\centering
\includegraphics[scale=1.2]{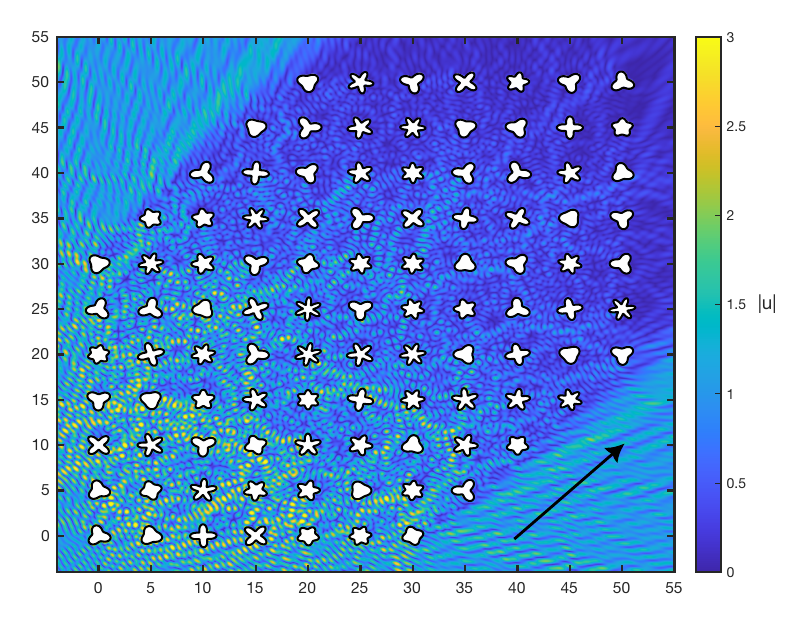}
\caption{ Flexural wave scattering by a collection of 101 starfish-shaped cavities in a thin elastic plate. Free plate boundary conditions are prescribed {on} each of the starfish, allowing the wave to pass quite far into the array before being attenuated. The wavenumber of the incident wave is $k = 6$ and the angle of the incident wave ($\theta = \pi/4$) is depicted by the black arrow.    }\label{multiplescattering}
\end{figure}

A convenient feature of the boundary integral equations described in this paper is that, after suitable discretization, they are amenable to standard techniques from fast algorithms such as {\it inter alia} fast multipole methods \cite{fmm1,fmm7}, and fast direct solvers~\cite{martinsson2019fast}. As a demonstration of this, we conclude our numerical illustrations with an example of the BIE applied to a large scale problem. In particular, we consider a plane wave scattering off of 101 {starfish-shaped} inclusions with free plate boundary  
conditions {(Figure~\ref{multiplescattering})}. The resulting discretized system has 76,480 unknowns. To accelerate the computation and reduce the overall memory cost we {compress low-rank interactions in the matrix using proxy surfaces, as in~\cite{martinsson2019fast}}. We remark that a `branching' structure is visible in the magnitude of the field. For simpler models this structure has been previously observed both experimentally and numerically in the context of flexural waves \cite{Jose2022, Jose2023, matula1995energy, Darabi2018}. Our formulation allows for further exploration of this phenomenon.

\clearpage
\section{Discussion} \label{discussion}

In this paper, we have presented integral equation methods for solving flexural wave problems for three common sets of boundary
conditions, either by extending existing methods as in the case of the clamped and supported plates, or by developing new methods in the case of the free plate. For the latter, we have used an integral representation that incorporates the composition of a standard layer potential with the Hilbert transform in order to obtain a second kind integral equation. Since the Hilbert transform offers another set of options for kernels in the integral representation, we believe this method may be useful more generally for the development of new integral representations for BVPs. Though not the focus of this study, interior problems can also be accurately solved using these integral representations. Given a compact geometry, one can
 study its resonant frequencies and resonant modes by computing the roots of the Fredholm determinant~\cite{ammari2009layer,zhao2015robust,Askham2020}. Further, the representations should be applicable to the static ($k=0$)
 case, a classical model in elasticity.

Our hope is that these methods will be useful to researchers interested in flexural waves
in both the applied engineering and glaciological communities. The methods are high-order accurate and 
can be readily extended to larger systems using existing fast algorithms. The evolution of 
sea ice and ice shelves is a multiscale phenomenon, and there is a growing desire {within this community to solve 
coupled problems quickly and} at multiple scales \cite{Banwell2023}. Fast algorithms based on second-kind integral equation representations are well-equipped to handle these large scattering problems, and we expect these methods to be effective in modeling wave propagation in {ice}. We anticipate this paper {to} be the first in a series that leverages integral equations to solve wave scattering problems related to the modeling of sea ice.

The boundary integral equations derived above were shown to be second kind for sufficiently smooth
geometries, requiring up to four absolutely continuous derivatives in the boundary parametrization
for the Taylor series estimates. There are a number of 
further questions to explore, including the behavior of the representations on geometries of lower
regularity; the solvability of the boundary integral equations and whether or not they are subject
to spurious resonances or spurious near-resonances~\cite{zhao2015robust}; and the uniqueness properties 
of solutions to the exterior free plate problem. The authors plan to explore the development of
representations for these problems with weaker dependence on the smoothness of the curve.
Such representations may be necessary for the treatment
of domains with geometric singularities, like corners and cracks. The solvability of the integral equations will
be treated in future work, including the characterization of any nullspaces and resonances of the equations
for both the $k\ne 0$ and $k=0$ cases. The application of these methods to transmitting flexural
waves is also being vigorously pursued.

\section{Code availability}

The integral equations were implemented and solved in \texttt{ChunkIE}, a MATLAB integral equations toolbox (\url{https://github.com/fastalgorithms/chunkie}). The examples from this paper may be found at \url{https://github.com/askhamwhat/flex-paper-examples}. 

\section{Acknowledgements}

We would like to thank {the two anonymous reviewers,} Douglas MacAyeal, Shidong Jiang, Tristan Goodwill, Mary Silber, {and Manas Rachh} for many useful {comments and} discussions. P.N. would like to acknowledge the support of the National Science Foundation (NSF) under Grant No. 2332479. J.G. Hoskins would like to thank the American Institute of Mathematics and, in particular, John Fry for hosting him on Bock Cay during the SQuaREs program, where parts of this work were completed. 

\appendix

\section{Heuristic strategy for deriving the integral kernels}\label{app:heuristic}

    Here we describe a strategy for deriving suitable integral kernels $K_i$ and
    describe its extension to the free and supported plate problems. In order to derive the kernels, it is useful to look at the Green's function on a domain with infinitely many degrees of symmetry. Consider the upper half-plane $\Omega := \{ (x_1,x_2) | x_2 \geq 0 \}$. In this geometry, the convolutions
    that define layer potentials can be understood in terms of the Fourier transforms of the 
    integral kernels in the $x_1$ direction; for an application, see~\cite{oneil2014efficient}. 
    The Green's function has the Sommerfeld integral identity~\cite{sommerfeld1949partial}:  
    \begin{align}   G(\bx ,\by ) &= {\frac{1}{2\pi}} \int^{\infty}_{-\infty      }  \hat{G}(\xi,x_2) \exp{(i\xi (x_1 - y_1 ))} \, \dd \xi \, , \label{sommerfeld} 
\end{align}
where
\begin{align}
    \hat{G}(\xi,x_2) &= {\frac{1}{4 k^2}}  \left ( \frac{\exp{(- x_2 \sqrt{\xi^2 - k^2}) }}{\sqrt{\xi^2 - k^2}} - \frac{\exp{(-x_2 \sqrt{\xi^2 + k^2 }) }}{\sqrt{\xi^2 + k^2}} \right ) \, ,
\end{align}
given $\bx  = (x_1,x_2)$ and $\by =(y_1,0)$. Because the expression for $\hat{G}$ is separated in the $x_1$ and $x_2$ directions,
it is simple to compute a similar Fourier transform of normal and tangential derivatives of the Green's function.

For any candidate kernel, $K_i$, given as some linear combination of derivatives of the Green's function, it is then possible 
to evaluate the $x_2\to 0^+$ limit of the boundary conditions applied to that kernel in terms of its 
Fourier transform.
Generally, the asymptotic expansion of this limit will include terms of the form $\xi^m$ or $\xi^m \sign(\xi)$ for $m \in \mathbb{Z}$ as 
$\xi \to  \pm \infty$ which can be used to characterize the corresponding boundary integral operator: 
positive values of $m$ correspond to hyper-singular boundary integral kernels in the $\mathbf{K}$ matrix or differential operators
in the $\mathbf{D}$ matrix and a non-zero constant term in the expansion corresponds to a constant 
jump while a signed constant corresponds to a Hilbert transform. The goal in integral kernel design is that this limiting
system should be second kind. Ideally, $D_{ii}+K_{ii}$ should have a constant
term and no higher order terms in the asymptotic expansion while $D_{ij}+K_{ij}$
for $i\ne j$ should have only $o(1)$ terms. 

For example, the $\xi \to \pm \infty$ asymptotic expansions of some derivatives of
$G$ are below 
\begin{align}
    \lim_{x_2\to 0^+} \widehat{G_{\bn_\by \bn_\by \bn_\by}} &= \frac{1}{2} + o(1) \, , \label{fourier1} \\
    \lim_{x_2\to 0^+} \widehat{G_{\bn_\by \btau_\by \btau_\by}} &= 0 + o(1) \, , \\
    \lim_{x_2\to 0^+} \widehat{G_{\bn_\bx \bn_\by \bn_\by \bn_\by}} &= \frac{3}{4} |\xi| + o(1)  \, , \\
    \lim_{x_2\to 0^+} \widehat{G_{\bn_\bx \bn_\by \btau_\by \btau_\by}} &= -\frac{1}{4} |\xi| + o(1) \, .     
\end{align}
From these, we see that the $K_1$ kernel for the clamped plate problem, 
$K_1 = G_{\bn_\by \bn_\by \bn_\by} + 3 G_{\bn_\by \btau_\by \btau_\by}$, should have 
\begin{align}
    \widehat{D_{11} + K_{11}} &= \frac{1}{2} + o(1) \, , \\
    \widehat{D_{21} + K_{21}} &= o(1) \, ,
\end{align}
as desired. 

Here we make a few observations. The first is that this analysis predicts,
correctly, that $D_{11} + \mathcal{K}_{11}$ should be second kind, with the constant in the
Fourier transform corresponding to the interior jump for this geometry. The second
is that the curvature-dependent term in $D_{21} + \mathcal{K}_{21}$ is not predicted because the
boundary is flat. Finally, we note that it is a general fact that the asymptotic
expansions contain terms of the form {$|\xi|^m$ if the number of tangential derivatives is
even and of the form $|\xi|^m\sign (\xi)$ if the number is odd. That is why the linear
combination of a term with zero tangential derivatives and a term with two tangential
derivatives is able to achieve cancellations in the singular parts.} 

The effect of curvature can be better understood by considering the analogous
problem on the disk, where the Fourier series of the Green's function can be written using Graf's addition formula: 
\begin{align}
    G(\bx , \by ) = \frac{1}{2k^2} \sum_{n=-\infty}^\infty \left[ \frac{i}{4} H_n^{(1)}(kr_\bx )J_n(kr_\by ) - \frac{1}{2\pi} K_n(kr_\bx ) I_n(kr_\by ) \right] \exp{(in\alpha)} \, , \label{grafs} 
\end{align}
subject to $r_\bx  > r_\by $, where $r_\bx $ is the target radius, $r_\by $ is the source radius, and $\alpha$ is the angle between the source and the target. As above, the separated representation allows for relatively
straightforward calculations of normal and tangential derivatives. The asymptotic expansions 
(as $n \to \pm \infty$) of the Fourier coefficients can then be computed and a similar analysis 
applied. If the disk radius, $r_\by$, appears in the asymptotics, this corresponds to a curvature
dependence. This can help guide the derivation of the jump relations on general curves, {though any terms involving the arc length derivative of curvature will fail to appear}.

For the free plate problem, we use a slightly different representation of the form
$\mathcal{K}_1 = \mathcal{K}_1^a + \mathcal{K}_1^b \mathcal{H}$. Because the Fourier
symbol of the Hilbert transform is $-i\sign(\xi)$, such a representation allows us to
cancel the singular parts of $K_1^a$ and $K_1^b$, even though $K_1^a$ has {no tangential
derivatives} and $K_1^b$ has {one}. As observed in Section~\ref{freesection}, this
effect can be understood for general curves by applying the Poincar\'{e}-Bertrand formula
and adding and subtracting the Hilbert transform from appropriate parts of the kernel. 
This Hilbert transform-based strategy then gives new degrees of freedom
in the design of integral kernels on general curves and surfaces and is likely applicable to 
integral kernel design for a variety of problems.

{
\section{Asymptotics of the Green's function}\label{app:asymptotics}

Our aim is to expand the flexural wave Green's function $G(\bx,\by)$ in the limit as $\|\bx - \by \| \to  0$ in order to understand its behavior near the diagonal. Recall that the Green's function is given by $G(\bx ,\by ) =  \displaystyle\frac{1}{2k^2} \bigg[\frac{i}{4} H_0^{(1)} (k \|\bx  - \by \|) - \frac{1}{2\pi} K_0 (k \| \bx  - \by \|) \bigg]$. 
Using the power series formulae for the Bessel functions $J_0$ and $Y_0$~\cite[\S 10.8]{NIST}, we have that 
\begin{equation}
H_0^{(1)} (k \|\bx  - \by \|) = \Psi_k(\bx ,\by ) + \frac{2i}{\pi} \ln \|\bx -\by \| \sum_{m\geq 0} 
\frac{(-1)^m}{(m!)^2} \left(\frac{k \|\bx -\by \|}{2}\right)^{2m} \, ,
\end{equation}
where $\Psi_k(\bx ,\by )$ is a smooth function. Likewise, 
\begin{equation}
 K_0 (k \|\bx  - \by \|) = \frac{\pi i}{2} \Psi_{ik}(\bx ,\by ) - \ln\|\bx -\by \| \sum_{m\geq 0} 
 \frac{1}{(m!)^2} \left(\frac{k\|\bx -\by \|}{2}\right)^{2m} \, .
\end{equation}
Combining these, we obtain  
\begin{align}
 G(\bx ,\by ) =  \Phi_{k}(\bx ,\by ) + \frac{1}{4\pi k^2}\ln\|\bx -\by \| \sum_{m\geq 0} 
 \frac{1-(-1)^m}{(m!)^2} \left(\frac{k\|\bx -\by \|}{2}\right)^{2m} \, ,
\end{align}
where $\Phi_k(\bx ,\by ) = \frac{i}{8k^2}( \Psi_k(\bx ,\by ) - \Psi_{ik}(\bx ,\by ))$ is
a smooth function. Separating out the first non-zero term from the sum, we have 
\begin{align}
  \label{eq:gseries}
 G(\bx ,\by ) =  \Phi_{k}(\bx ,\by ) + \frac{1}{8\pi} \|\bx -\by \|^2 \ln\|\bx -\by \| + \frac{1}{4\pi k^2}\ln\|\bx -\by \| \sum_{m\geq 3} 
 \frac{1-(-1)^m}{(m!)^2} \left(\frac{k\|\bx -\by \|}{2}\right)^{2m} \, ,
\end{align}
so that the leading non-smooth term is precisely the Green's function of the biharmonic equation. Additionally, the next order term takes the form $\|\bx -\by \|^6 \ln\|\bx -\by \|$. Since we are taking at most five derivatives in the present study, we can restrict our analysis of singular behavior to the biharmonic Green's function since the rest of the terms in the asymptotic expansions will be continuous. 
} 

\section{Formulae for the derivatives of the biharmonic Green's function} 
In this section, we denote the biharmonic Green's function by $G^{B}(\bx,\by) = \frac{1}{16\pi} \|\br \|^2 \ln\|\br\|^2 $, where $\br = \bx - \by$. For completeness, we provide formulae for the derivatives of the biharmonic Green's function, which are used to analyze the leading order behavior of the flexural wave kernels. {We provide the corresponding kernels for the clamped plate problem in \ref{appclamped}, the supported plate problem in \ref{appsupported}, and the free plate problem in \ref{appfree}.}

\subsection{Clamped plate kernels} \label{appclamped}
First, we provide the following derivatives of $G^{B}$ used in the integral equation for the clamped plate (Theorem \ref{lem:clampedsmooth}). A straightforward but tedious calculation shows,
\begin{align}
    G^{B}_{\bn_\by  \bn_\by }{(\bx,\by)} &= \frac{1}{4\pi} \frac{[\br\cdot  \bn (\by) ]^2}{\|\br \|^2} + \frac{1}{8\pi}\ln
    \|\br \|^2 + \frac{1}{8\pi} \, , \\
    G^{B}_{\btau_\by  \btau_\by } {(\bx,\by)} &= \frac{1}{4\pi} \frac{[\br\cdot \btau (\by)]^2}{\|\br \|^2} + \frac{1}{8\pi}\ln\|\br \|^2 + \frac{1}{8\pi} \, ,\\
    G^{B}_{\bn_\bx  \bn_\by  \bn_\by } {(\bx,\by)} &= \frac{1}{2\pi} \frac{[\br\cdot \bn (\by)][\bn (\by) \cdot  \bn (\bx ) ]}{\|\br \|^2} - \frac{1}{2\pi} \frac{[\br\cdot  \bn (\bx ) ][\br\cdot \bn (\by)]^2}{\|\br \|^4} +\frac{1}{4\pi}\frac{\br\cdot  \bn (\bx ) }{\|\br \|^2} \, ,\\
    G^{B}_{\bn_\bx  \btau_\by  \btau_\by } {(\bx,\by)} &= \frac{1}{2\pi} \frac{[\br\cdot \btau (\by)][\btau (\by) \cdot  \bn (\bx ) ]}{\|\br \|^2} - \frac{1}{2\pi} \frac{[\br\cdot  \bn (\bx ) ][\br\cdot \btau (\by)]^2}{\|\br \|^4} +\frac{1}{4\pi}\frac{\br\cdot  \bn (\bx ) }{\|\br \|^2} \, ,\\
    G^{B}_{\bn_\by  \bn_\by  \bn_\by } {(\bx,\by)} &= -\frac{3}{4\pi} \frac{[\br\cdot \bn (\by)]}{\|\br \|^2} + \frac{1}{2\pi} \frac{[\br\cdot \bn (\by)]^3}{\|\br \|^4} \, , \\
    G^{B}_{ \bn_\by  \btau_\by  \btau_\by } {(\bx,\by)} &=  \frac{1}{2\pi} \frac{[\br\cdot  \bn (\by) ][\br\cdot \btau (\by)]^2}{\|\br \|^4} -\frac{1}{4\pi}\frac{\br\cdot  \bn (\by) }{\|\br \|^2} \ ,\\
      G^{B}_{\bn_\bx  \bn_\by  \bn_\by  \bn_\by } {(\bx,\by)} &= -\frac{3}{4\pi} \frac{ \bn (\bx ) \cdot  \bn (\by) }{\|\br \|^2} + \frac{3}{2\pi} \frac{[\br\cdot  \bn (\by) ][\br\cdot  \bn (\bx ) ]}{\|\br \|^4}+ \frac{3}{2\pi} \frac{[\br\cdot  \bn (\by) ]^2 [ \bn (\bx )  \cdot  \bn (\by) ]}{\|\br \|^4} \nonumber \\
      &\qquad -\frac{2}{\pi} \frac{[\br \cdot \bn (\by)]^3[\br\cdot \bn (\bx )]}{\|\br \|^6} \, , \\
    G^B_{\bn_\bx  \bn_\by  \btau_\by  \btau_\by } {(\bx,\by)} &=  \frac{1}{2\pi} \frac{[\bn (\bx ) \cdot \bn (\by)][\br\cdot \btau (\by)]^2}{\|\br \|^4} + \frac{1}{\pi} \frac{[\br\cdot \bn (\by)][\br\cdot \btau (\by)][\bn (\bx ) \cdot \btau (\by)]}{\|\br \|^4} \nonumber \\
    &\qquad - \frac{2}{\pi} \frac{[\br\cdot \bn (\by)][\br\cdot \btau (\by)]^2 [\br\cdot \bn (\bx )]}{\|\br \|^6} - \frac{1}{4\pi} \frac{\bn (\bx ) \cdot \bn (\by) }{\|\br \|^2} \nonumber \\
    &\qquad + \frac{1}{2\pi} \frac{[\br\cdot \bn (\by)][\br\cdot \bn (\bx )]}{\|\br \|^4}  \, .
\end{align}

\pagebreak

\subsection{Supported plate kernels} \label{appsupported}

{Next}, we provide the following derivatives of $G^{B}$ used in the integral equation for the supported plate (Theorem \ref{supportedsmooth}),
\begin{align}
        G^{B}_{\bn_\by  } {(\bx,\by)} &= - \frac{1}{8\pi} [\br \cdot \bn (\by)] \ln\|\br \|^2 -  \frac{1}{8\pi} [\br \cdot \bn (\by)]  \, ,  \\
        G^{B}_{\btau_\by  } {(\bx,\by)} &= - \frac{1}{8\pi} [\br \cdot \btau (\by)] \ln\|\br \|^2 -  \frac{1}{8\pi} [\br \cdot \btau (\by)] \, ,   \\
        G^{B}_{\bn_\by  \bn_\by } {(\bx,\by)} &= \frac{1}{4\pi} \frac{[\br\cdot  \bn (\by) ]^2}{\|\br \|^2} + \frac{1}{8\pi}\ln\|\br \|^2 + \frac{1}{8\pi} \, , \\
    G^{B}_{\bn_\by  \bn_\by  \bn_\by } {(\bx,\by)} &= -\frac{3}{4\pi} \frac{[\br\cdot \bn (\by)]}{\|\br \|^2} + \frac{1}{2\pi} \frac{[\br\cdot \bn (\by)]^3}{\|\br \|^4} \, ,\\
    G^{B}_{ \bn_\by  \btau_\by  \btau_\by } {(\bx,\by)} &=  \frac{1}{2\pi} \frac{[\br\cdot  \bn (\by) ][\br\cdot \btau (\by)]^2}{\|\br \|^4} -\frac{1}{4\pi}\frac{\br\cdot  \bn (\by) }{\|\br \|^2} \, , \\
  G^{B}_{\bn_\bx  \bn_\bx  \bn_\by } {(\bx,\by)} &= -\frac{1}{2\pi} \frac{[\br\cdot \bn (\bx )][\bn (\bx ) \cdot  \bn (\by) ]}{\|\br \|^2} + \frac{1}{2\pi} \frac{[\br\cdot  \bn (\by) ][\br\cdot \bn (\bx )]^2}{\|\br \|^4} -\frac{1}{4\pi}\frac{\br\cdot  \bn (\by) }{\|\br \|^2} \, , \\
    G^{B}_{\btau_\bx  \btau_\bx  \bn_\by } {(\bx,\by)} &= -\frac{1}{2\pi} \frac{[\br\cdot \btau (\bx)][\btau (\bx) \cdot  \bn (\by) ]}{\|\br \|^2} + \frac{1}{2\pi} \frac{[\br\cdot  \bn (\by) ][\br\cdot \btau (\bx)]^2}{\|\br \|^4} -\frac{1}{4\pi}\frac{\br\cdot  \bn (\by) }{\|\br \|^2} \, , \\
    G^{B}_{\bn_\bx  \bn_\bx  \btau_\by } {(\bx,\by)} &= -\frac{1}{2\pi} \frac{[\br\cdot  \bn (\bx ) ][ \bn (\bx )  \cdot  \btau (\by) ]}{\|\br \|^2} + \frac{1}{2\pi} \frac{[\br\cdot  \btau (\by) ][\br\cdot  \bn (\bx ) ]^2}{\|\br \|^4} -\frac{1}{4\pi}\frac{\br\cdot  \btau (\by) }{\|\br \|^2} \, ,\\
     G^{B}_{\btau_\bx  \btau_\bx  \btau_\by } {(\bx,\by)} &= -\frac{1}{2\pi} \frac{[\br\cdot \btau (\bx)][\btau (\bx) \cdot  \btau (\by) ]}{\|\br \|^2} + \frac{1}{2\pi} \frac{[\br\cdot  \btau (\by) ][\br\cdot \btau (\bx)]^2}{\|\br \|^4} -\frac{1}{4\pi}\frac{\br\cdot  \btau (\by) }{\|\br \|^2}\, ,\\
    G^{B}_{\bn_\bx  \bn_\bx  \bn_\by  \bn_\by } {(\bx,\by)} &= \frac{1}{2\pi} \frac{[\bn (\bx )\cdot \bn (\by)]^2}{\|\br \|^2} - \frac{2}{\pi} \frac{[\br\cdot \bn (\bx )][\br \cdot \bn (\by)][\bn (\bx )\cdot \bn (\by)]}{\|\br \|^4} + \frac{1}{4\pi} \frac{1}{\|\br \|^2} \nonumber \\
    &\qquad -\frac{1}{2\pi}\frac{[\br\cdot \bn (\by)]^2}{\|\br \|^4} + \frac{2}{\pi} \frac{[\br\cdot \bn (\bx )]^2[\br\cdot \bn (\by)]^2}{\|\br \|^6} - \frac{1}{2\pi} \frac{[\br\cdot \bn (\bx )]^2}{\|\br \|^4} \, , \\
    G^{B}_{\btau_\bx  \btau_\bx  \bn_\by  \bn_\by } {(\bx,\by)} &= \frac{1}{2\pi} \frac{[\btau (\bx)\cdot \bn (\by)]^2}{\|\br \|^2} - \frac{2}{\pi} \frac{[\br\cdot \btau (\bx)][\br \cdot \bn (\by)][\btau (\bx)\cdot \bn (\by)]}{\|\br \|^4} + \frac{1}{4\pi} \frac{1}{\|\br \|^2} \nonumber \\
    &\qquad -\frac{1}{2\pi}\frac{[\br\cdot \bn (\by)]^2}{\|\br \|^4} + \frac{2}{\pi} \frac{[\br\cdot \btau (\bx)]^2[\br\cdot \bn (\by)]^2}{\|\br \|^6} - \frac{1}{2\pi} \frac{[\br\cdot \btau (\bx)]^2}{\|\br \|^4}  \, , 
\end{align}
\begin{align}
  G^B_{\bn_{\bx } \bn_{\bx } \bn_{\by } \bn_{\by } \bn_{\by }} {(\bx,\by)} &= \frac{3}{2 \pi} \frac{\br \cdot \bn (\by)}{\|\br \|^4} - \frac{6}{\pi} \frac{[\br \cdot \bn (\by)][\br\cdot \bn (\bx )]^2}{\|\br \|^6} + \frac{3}{\pi} \frac{[\br \cdot \bn (\by)][\bn (\bx ) \cdot \bn (\by)]^2}{\|\br \|^4} \nonumber \\
    &\qquad - \frac{12}{\pi} \frac{[\br\cdot \bn (\by)]^2 [\br \cdot \bn( \bx ) ][\bn (\bx ) \cdot \bn (\by)]}{\|\br \|^6} - \frac{2}{\pi} \frac{[\br \cdot \bn (\by)]^3}{\|\br \|^6} \nonumber \\
    &\qquad + \frac{12}{\pi} \frac{[\br \cdot \bn (\by)]^3 [\br\cdot \bn (\bx )]^2}{\|\br \|^8} + \frac{3}{\pi} \frac{[\br\cdot \bn (\bx )][\bn (\bx ) \cdot \bn (\by)]}{\|\br \|^4}\, , \label{Gnxnxnynyny}  \\
    G^B_{\btau_{\bx } \btau_{\bx } \bn_{\by } \bn_{\by } \bn_{\by }} {(\bx,\by)} &= \frac{3}{2 \pi} \frac{\br \cdot \bn (\by)}{\|\br \|^4} - \frac{6}{\pi} \frac{[\br \cdot \bn (\by)][\br\cdot \btau (\bx)]^2}{\|\br \|^6} + \frac{3}{\pi} \frac{[\br \cdot \bn (\by)][\btau (\bx) \cdot \bn (\by)]^2}{\|\br \|^4} \nonumber \\
    &\qquad - \frac{12}{\pi} \frac{[\br\cdot \bn (\by)]^2 [\br \cdot \btau( \bx ) ][\btau (\bx) \cdot \bn (\by)]}{\|\br \|^6} - \frac{2}{\pi} \frac{[\br \cdot \bn (\by)]^3}{\|\br \|^6} \nonumber \\
    &\qquad + \frac{12}{\pi} \frac{[\br \cdot \bn (\by)]^3 [\br\cdot \btau (\bx)]^2}{\|\br \|^8}+ \frac{3}{\pi} \frac{[\br\cdot \btau (\bx)][\btau (\bx) \cdot \bn (\by)]}{\|\br \|^4} \, , \\
    G^B_{\bn_\bx  \bn_\bx  \bn_\by  \btau_\by  \btau_\by } {(\bx,\by)} &= \frac{2}{\pi} \frac{[\bn (\bx ) \cdot \bn (\by) ] [\bn (\bx ) \cdot \btau (\by) ] [\br \cdot \btau (\by)]}{\|\br \|^4} - \frac{4}{\pi} \frac{[\bn (\bx ) \cdot \bn (\by)][\br\cdot \btau (\by)]^2 [\br\cdot \bn (\bx )]}{\|\br \|^6} \nonumber \\
    &\qquad + \frac{1}{\pi} \frac{[\br \cdot \bn (\by) ] [\bn (\bx ) \cdot \btau (\by)]^2}{\|\br \|^4} - \frac{8}{\pi} \frac{[\br \cdot \bn (\by)] [\br \cdot \btau (\by)] [\br\cdot \bn (\bx )] [\bn (\bx ) \cdot \btau (\by)]}{\|\br \|^6} \nonumber \\
    &\qquad -\frac{2 }{\pi} \frac{[\br \cdot \bn (\by)] [\br\cdot \btau (\by)]^2}{\|\br \|^6} + \frac{12}{\pi} \frac{[\br \cdot \bn (\by)][\br\cdot \btau (\by)]^2 [\br\cdot \bn (\bx )]^2}{\|\br \|^8} \nonumber \\
    &\qquad +\frac{1}{\pi} \frac{[\bn (\bx ) \cdot \bn (\by)] [\br\cdot \bn (\bx )]}{\|\br \|^4} + \frac{1}{2\pi} \frac{\br \cdot \bn (\by)}{\|\br \|^4}  - \frac{2}{\pi} \frac{[\br\cdot \bn (\by)][\br\cdot \bn (\bx )]^2}{\|\br \|^6} \, , \\
    G^B_{\btau_\bx  \btau_\bx  \bn_\by  \btau_\by  \btau_\by } {(\bx,\by)} &= \frac{2}{\pi} \frac{[\btau (\bx) \cdot \bn (\by) ] [\btau (\bx) \cdot \btau (\by) ] [\br \cdot \btau (\by)]}{\|\br \|^4} - \frac{4}{\pi} \frac{[\btau (\bx) \cdot \bn (\by)][\br\cdot \btau (\by)]^2 [\br\cdot \btau (\bx)]}{\|\br \|^6} \nonumber \\
    &\qquad + \frac{1}{\pi} \frac{[\br \cdot \bn (\by) ] [\btau (\bx) \cdot \btau (\by)]^2}{\|\br \|^4} - \frac{8}{\pi} \frac{[\br \cdot \bn (\by)] [\br \cdot \btau (\by)] [\br\cdot \btau (\bx)] [\btau (\bx) \cdot \btau (\by)]}{\|\br \|^6} \nonumber \\
    &\qquad -\frac{2 }{\pi} \frac{[\br \cdot \bn (\by)] [\br\cdot \btau (\by)]^2}{\|\br \|^6} + \frac{12}{\pi} \frac{[\br \cdot \bn (\by)][\br\cdot \btau (\by)]^2 [\br\cdot \btau (\bx)]^2}{\|\br \|^8} \nonumber \\
    &\qquad +\frac{1}{\pi} \frac{[\btau (\bx) \cdot \bn (\by)] [\br\cdot \btau (\bx)]}{\|\br \|^4} + \frac{1}{2\pi} \frac{\br \cdot \bn (\by)}{\|\br \|^4}  - \frac{2}{\pi} \frac{[\br\cdot \bn (\by)][\br\cdot \btau (\bx)]^2}{\|\br \|^6} \, .
\end{align}

\pagebreak

\subsection{Free plate kernels} \label{appfree}
{Lastly}, we provide the following derivatives of $G^{B}$ used in the integral equation for the clamped plate (Theorem \ref{freesmooth}). Again, a straightforward but tedious calculation shows,
\begin{align}
    G^{B}_{\bn_\bx  \bn_\bx } {(\bx,\by)} &= \frac{1}{4\pi} \frac{[\br\cdot  \bn (\bx ) ]^2}{\|\br \|^2} + \frac{1}{8\pi}\ln\|\br \|^2 + \frac{1}{8\pi} \, , \\
    G^{B}_{\btau_\bx  \btau_\bx } {(\bx,\by)} &= \frac{1}{4\pi} \frac{[\br\cdot \btau (\bx)]^2}{\|\br \|^2} + \frac{1}{8\pi}\ln\|\br \|^2 + \frac{1}{8\pi} \, , \\
    G^{B}_{\bn_\bx  \bn_\bx  \bn_\by } {(\bx,\by)} &= -\frac{1}{2\pi} \frac{[\br\cdot \bn (\bx )][\bn (\bx ) \cdot  \bn (\by) ]}{\|\br \|^2} + \frac{1}{2\pi} \frac{[\br\cdot  \bn (\by) ][\br\cdot \bn (\bx )]^2}{\|\br \|^4} -\frac{1}{4\pi}\frac{\br\cdot  \bn (\by) }{\|\br \|^2} \, ,\\
    G^{B}_{\btau_\bx  \btau_\bx  \bn_\by } {(\bx,\by)} &= -\frac{1}{2\pi} \frac{[\br\cdot \btau (\bx)][\btau (\bx) \cdot  \bn (\by) ]}{\|\br \|^2} + \frac{1}{2\pi} \frac{[\br\cdot  \bn (\by) ][\br\cdot \btau (\bx)]^2}{\|\br \|^4} -\frac{1}{4\pi}\frac{\br\cdot  \bn (\by) }{\|\br \|^2} \, , \\
    G^{B}_{\bn_\bx  \bn_\bx  \btau_\by } {(\bx,\by)} &= -\frac{1}{2\pi} \frac{[\br\cdot  \bn (\bx ) ][ \bn (\bx )  \cdot  \btau (\by) ]}{\|\br \|^2} + \frac{1}{2\pi} \frac{[\br\cdot  \btau (\by) ][\br\cdot  \bn (\bx ) ]^2}{\|\br \|^4} -\frac{1}{4\pi}\frac{\br\cdot  \btau (\by) }{\|\br \|^2} \, , \\
     G^{B}_{\btau_\bx  \btau_\bx  \btau_\by } {(\bx,\by)} &= -\frac{1}{2\pi} \frac{[\br\cdot \btau (\bx)][\btau (\bx) \cdot  \btau (\by) ]}{\|\br \|^2} + \frac{1}{2\pi} \frac{[\br\cdot  \btau (\by) ][\br\cdot \btau (\bx)]^2}{\|\br \|^4} -\frac{1}{4\pi}\frac{\br\cdot  \btau (\by) }{\|\br \|^2} \, , \\
    G^{B}_{\bn_\bx  \bn_\bx  \bn_\bx } {(\bx,\by)} &= \frac{3}{4\pi} \frac{\br\cdot \bn (\bx )}{\|\br \|^2} - \frac{1}{2\pi} \frac{[\br\cdot \bn (\bx )]^3}{\|\br \|^4} \, , \\
    G^{B}_{\bn_\bx  \btau_\bx  \btau_\bx  } {(\bx,\by)} &=  - \frac{1}{2\pi} \frac{[\br\cdot  \bn (\bx ) ][\br\cdot \btau (\bx)]^2}{\|\br \|^4} +\frac{1}{4\pi}\frac{\br\cdot  \bn (\bx ) }{\|\br \|^2} \, , \\
    G^{B}_{\bn_\bx  \bn_\bx  \bn_\bx  \bn_\by } {(\bx,\by)} &= -\frac{3}{4\pi} \frac{ \bn (\bx ) \cdot  \bn (\by) }{\|\br \|^2} + \frac{3}{2\pi} \frac{[\br\cdot  \bn (\by) ][\br\cdot  \bn (\bx ) ]}{\|\br \|^4}+ \frac{3}{2\pi} \frac{[\br\cdot  \bn (\bx ) ]^2 [ \bn (\bx )  \cdot  \bn (\by) ]}{\|\br \|^4} \nonumber \\
    &\qquad -\frac{2}{\pi} \frac{[\br \cdot \bn (\bx )]^3[\br\cdot \bn (\by)]}{\|\br \|^6} \, , \\
    G^{B}_{\bn_\bx  \btau_\bx  \btau_\bx   \bn_\by } {(\bx,\by)} &= \frac{1}{2\pi} \frac{[\br \cdot \btau (\bx)]^2 [ \bn (\bx )  \cdot  \bn (\by) ]}{\|\br \|^4} + \frac{1}{\pi} \frac{[\br\cdot \btau (\bx) ][\btau (\bx)\cdot  \bn (\by) ] [\br \cdot  \bn (\bx ) ] }{\|\br \|^4} \nonumber \\
    &\qquad - \frac{2}{\pi} \frac{[\br \cdot  \bn (\by) ][\br \cdot  \bn (\bx ) ] [\br \cdot \btau (\bx)]^2 }{\|\br \|^6} -\frac{1}{4\pi} \frac{ \bn (\bx )  \cdot  \bn (\by) }{\|\br \|^2 } \nonumber \\ &\qquad+ \frac{1}{2\pi} \frac{[\br \cdot  \bn (\by) ][\br \cdot  \bn (\bx ) ]}{\|\br \|^4} \, ,
\end{align}
\begin{align}
     G^{B}_{\bn_\bx  \bn_\bx  \bn_\bx  \btau_\by } {(\bx,\by)} &= -\frac{3}{4\pi}\frac{ \bn (\bx )  \cdot  \btau (\by) }{\|\br \|^2} + \frac{3}{2\pi} \frac{[\br\cdot  \bn (\bx ) ][\br\cdot  \btau (\by) ]}{\|\br \|^4} + \frac{3}{2\pi} \frac{[\br\cdot  \bn (\bx ) ]^2[  \bn (\bx )  \cdot  \btau (\by) ]}{\|\br \|^4} \nonumber \\
    &\qquad -\frac{2}{\pi}\frac{[\br\cdot  \btau (\by)  ][\br \cdot  \bn (\bx ) ]^3}{\|\br \|^6} \, ,\\
     G^{B}_{\bn_\bx  \btau_\bx  \btau_\bx  \btau_\by  } {(\bx,\by)} &= \frac{1}{\pi} \frac{[\br\cdot  \bn (\bx )  ][\br\cdot \btau (\bx)] [\btau (\bx)\cdot  \btau (\by) ]}{\|\br \|^4}+ \frac{1}{2\pi }\frac{[ \btau (\by)  \cdot  \bn (\bx ) ] [\br \cdot \btau (\bx)]^2}{\|\br \|^4} \nonumber \\
    &\qquad - \frac{2}{\pi} \frac{[\br\cdot  \bn (\bx ) ][\br \cdot  \btau (\by) ][\br \cdot \btau (\bx)]^2}{\|\br \|^6} - \frac{1}{4\pi} \frac{ \btau (\by)  \cdot  \bn (\bx ) }{\|\br \|^2} \nonumber \\
    &\qquad + \frac{1}{2\pi} \frac{[\br \cdot  \bn (\bx ) ] [\br \cdot  \btau (\by)  ]}{\|\br \|^4}  \, .
\end{align}

{
\section{Derivation of jump relations of the layer potentials}

The aim of this section is to establish the limit of several key boundary integral operators as the target $\bx$ approaches the boundary $\partial \Omega$. These limits, known as jump relations, are described in Theorems \ref{thm:clamped}, \ref{thm:supported}, and \ref{thm:free} for the clamped, supported, and free plates, respectively. In general, the integral operators can be split into a singular part and a principle value part, i.e. 
\begin{align}
    \lim_{\bx \to \bx_0^\pm}\mathcal{K}[\sigma](\bx) &=  \lim_{\bx \to \bx_0^\pm} \int_{\Gamma_\delta } K (\bx,\by) \sigma(\by) \, \dd S(\by) + \int_{\partial\Omega \setminus \Gamma_\delta} K (\bx_0,\by) \sigma(\by) \, \dd S(\by) \nonumber \\
    &= \lim_{\delta \to 0} \lim_{\bx \to \bx_0^\pm} \int_{\Gamma_\delta } K (\bx,\by) \sigma(\by) \, \dd S(\by) + \pv \int_{\partial\Omega} K (\bx_0,\by) \sigma(\by) \, \dd S(\by) \, ,
\end{align}
where $\bx_0 \in \partial\Omega$, $\bx_0^+$ and $\bx_0^-$ denote the exterior and interior limits respectively, and  $\Gamma_\delta = \left\{ \bgamma(s) \ | \ s \in [-\delta,\delta], \ \bgamma(0) = \bx_0 \right\}$ is a small region of $\partial \Omega$ containing $\bx_0 = \bgamma(0)$. 

In order to determine the contribution of the singular part, we consider an off-surface target $\bx_0 + h \bn(\bx_0) $, expand the kernels about $\bx_0$, and evaluate the limit as $h \to 0^\pm$, where $h > 0$ is the exterior limit and $h < 0$ is the interior limit. 
For $\bgamma$ sufficiently differentiable and $s$ sufficiently small, we can write:
\begin{align}
\bgamma(s) &= \bgamma(0) + s \btau(0) - \frac{s^2}{2} \kappa(0)\bn(0) -\frac{s^3}{6}(\kappa'(0) \bn(0)+\kappa^2(0) \btau(0)) + \mathcal{O}(s^4) \, , \\
\btau(s)  &= \btau(0) - s \kappa(0) \bn(0)- \frac{s^2}{2} (\kappa'(0) \bn(0) +\kappa^2(0) \btau(0)) + \mathcal{O}(s^3) \, , \\
\bn(s) &= \bn(0) + s \kappa(0) \btau(0)+ \frac{s^2}{2} ( \kappa'(0) \btau(0) - \kappa^2(0) \bn(0)) + \mathcal{O}(s^3) \, .   
\end{align}

\begin{lem}
For $h \in \mathbb{R}$, let $\bx = \bgamma(0) + h \bn(0)$ and $\by  = \bgamma(s)$. Then
\begin{align}
&\br = h \bn(0) -s \btau(0) +\frac{s^2}{2} \kappa(0) \bn(0) + \frac{s^3}{6} (\kappa'(0) \bn (0) + \kappa^2(0) \btau(0)) + \mathcal{O}(s^4) \, .
\end{align}
Additionally, we have the following expansions for terms that commonly appear in the kernels:
\begin{align}
&\|\br \|^2 = h^2 +s^2(1+h \kappa(0))+h\frac{s^3}{3} \kappa'(0) +\mathcal{O}(s^4) \, , \\
&\br \cdot \bn(0) = h + \frac{s^2}{2} \kappa(0) + \frac{s^3}{6} \kappa'(0) + \mathcal{O}(s^4) \, , \\
&\br \cdot \bn(s) = h - \frac{s^2}{2}(\kappa(0)+h \kappa^2(0)) + \mathcal{O}(s^3) \, , \\
&\br \cdot \btau(0) = -s+ \frac{\kappa^2(0)s^3}{6} + \mathcal{O}(s^4) \, , \\
&\br \cdot \btau(s) = -(1+h \kappa(0))s -\frac{s^2}{2} h \kappa'(0) + \mathcal{O}(s^3) \, , \\
&\bn(0) \cdot \bn(s) = \btau(0) \cdot \btau(s) = 1-\frac{\kappa^2(0)s^2}{2}  + \mathcal{O}(s^3) \, , \\
&\bn(0) \cdot \btau(s) = - \bn(s) \cdot \btau(0) = - \kappa(0) s - \kappa'(0) \frac{s^2}{2} + \mathcal{O}(s^3) \, ,
\end{align}
where the constant implicit in the big $\mathcal{O}$ notation above can be chosen independent of $h$ for $|h|<h_0.$
\end{lem}

Next, we define 
\begin{align}
R(s) &:= \|\br \|^2-h^2 \, .
\end{align}
We perform a standard ``blowing up of the singularity'', setting $s = uv(u)$ and $R(uv(u)) = u^2.$ Rearranging yields
\begin{align}
 (1+h\kappa(0)) v^2 + \frac13 h \kappa'(0) u v^3 + \mathcal{O}(u^2 v^4) = 1 \, . \label{blowup}
\end{align}
When $u = 0$, the equation is a quadratic with two distinct roots. Choosing the positive root, we seek to find a smooth function $v(u)$ satisfying 
\begin{align}
v(0)= \frac{1}{\sqrt{1+h\kappa(0)}} \, .
\end{align} 
It follows immediately from the implicit function theorem that there exists $h_0,u_0>0$ such that for all $|h|<h_0$ and all $|u|<u_0,$ $v(u,h)$ is well-defined and smooth. In particular, we can write $s(u) = v(u) u$ on the same neighborhood. The first few terms of this change of variables are as follows
\begin{align}
s(u) = \frac{1}{\sqrt{1+h\kappa(0)}} u  - \frac{h \kappa'(0) }{6 (1+h\kappa(0))^2} u^2 + \mathcal{O}(u^3) \, .
\end{align}
{It is easily shown that the constant implicit in the remainder term is bounded uniformly for all $|h|<h_0$ and that the remainder is smooth in $u$ with a derivative which is $\mathcal{O}(u^2)$ uniformly for $|h| <h_0.$} The coefficients of the terms above can be expanded in $h,$ yielding
\begin{align}
    s(u) = u-\frac{\kappa(0)}{2}uh+ \frac{3\kappa^2(0)}{8}uh^2  -\frac{\kappa'(0)}{6}u^2 h + \frac{\kappa(0)\kappa'(0)}{3}u^2h^2 + \mathcal{O}(u^3+h^3) \, , \label{s(u)}
\end{align}
and
\begin{align}
    s'(u) = 1 -\frac{\kappa(0)}{2}h-\frac{\kappa'(0)}{3} u h + \mathcal{O}(u^2+h^2) \, . \label{s'(u)}
\end{align}

\begin{lem}\label{lem:uexpansions} Using this change of variables, let $\by =  \bgamma(u)$. We now have the following expansions in terms of $u$:
\begin{align}
   & \|\br\|^2 = u^2 + h^2 \, , \\
   &  \br \cdot \bn(0) = h+ \frac{\kappa(0)}{2}u^2-\frac{\kappa^2(0)}{2}hu^2+\frac{\kappa^3(0)}{2}h^2u^2+\mathcal{O}(u^3+h^3) \, , \\
 &    \br \cdot \bn(s(u)) = h- \frac{\kappa(0)u^2}{2} +\mathcal{O}(u^3+ h^3) \, , \\
  &  \br \cdot \btau(0) = -u +\frac{\kappa(0)}{2}hu-\frac{3\kappa^2(0)}{8}h^2u+\frac{\kappa'(0)}{6}hu^2-\frac{\kappa(0)\kappa'(0)}{3}h^2u^2+ \mathcal{O}(u^3+h^3) \, , \\
  &      \br \cdot \btau(s(u)) = -u -\frac{\kappa(0)}{2}hu+\frac{\kappa^2(0)}{8}h^2u-\frac{\kappa'(0)}{3}hu^2-\frac{\kappa(0)\kappa'(0)}{3}h^2u^2+ \mathcal{O}(u^3+h^3) \, , \\
  &  \bn(0) \cdot \bn(s(u)) = \btau(0) \cdot \btau(s(u)) = 1-\frac{\kappa^2(0)}{2}u^2+\frac{\kappa^3(0)}{2}hu^2- \frac{\kappa^4(0)}{2}h^2 u^2 + \mathcal{O}(u^3+h^3) \, , \\
& \bn(0) \cdot \btau(s(u)) = - \bn(s(u)) \cdot \btau(0) = -\kappa(0) u +\frac{\kappa^2(0)}{2}hu-\frac{3\kappa^3(0)}{8}h^2 u - \frac{\kappa'(0)}{2}u^2 \nonumber \\
&\qquad\qquad\qquad\qquad+ \frac{2 \kappa(0)\kappa'(0)}{3} hu^2 - \frac{5 \kappa^2(0)\kappa'(0)}{6} h^2 u^2 + \mathcal{O}(u^3+h^3) \, .
\end{align}
\end{lem}

For the singular integrals involved in our integral representations, the change of variables takes the form
\begin{align}
\int^{\epsilon v(\epsilon)}_{-\epsilon v(\epsilon)} \frac{f(h,s)}{(h^2 + R(s))^n} \sigma(s) \, \dd s = \int^\epsilon_{-\epsilon} \frac{f(h,s(u))}{(h^2 + u^2)^n} \sigma(s(u)) s'(u) \, \dd u \, , \label{changeofvar}
\end{align}
for some $n \in \mathbb{N}.$ Here we choose $\epsilon > 0$ such that $[-\epsilon v(\epsilon),\epsilon v(\epsilon)] \subset [-\delta,\delta].$ It is straightforward to argue that the integral over $[-\delta,\delta] \setminus [-\epsilon v(\epsilon),\epsilon v(\epsilon)]$ is bounded uniformly for all $|h|<h_0$. 

\begin{lem}\label{lem:sigexpansion}
Assume that $\sigma$ is sufficiently smooth so that we can expand $\sigma(s)$ about $s = 0$:
\begin{align}
    \sigma(s) = \sigma(0) + s \sigma'(0) + \frac{s^2}{2} \sigma''(0) + \mathcal{O}(s^3) \, .
\end{align}
Then, using \eqref{s(u)} to perform the change of variables, we have:
\begin{align}
    \sigma(s(u)) &= \sigma(0) + u \sigma'(0) -\frac{\kappa(0)}{2}uh \sigma'(0) + \frac{3\kappa^2(0)}{8}uh^2 \sigma'(0)  -\frac{\kappa'(0)}{6}u^2 h \sigma'(0) \nonumber \\
    &\qquad+ \frac{\kappa(0)\kappa'(0)}{3}u^2h^2 \sigma'(0) + \frac{1}{2} u^2 \sigma''(0) - \frac{1}{2} h \kappa(0) u^2 \sigma''(0) + \frac{1}{2} h^2 u^2 \kappa^2(0) \sigma''(0) + \mathcal{O}(u^3 + h^3) \, .    
\end{align}
\end{lem}

The integrand in \eqref{changeofvar} may take on a few different forms. The next remark gives the list of integrals that arise from these expansions.
\begin{rmk} \label{rmk:integrals}
The following integrals will be used below:
\begin{align} 
&\int_{-\epsilon}^\epsilon \frac{h}{h^2+u^2}\, \dd u = 2 \arctan \left(\frac{\epsilon}{h}\right) \, , \\
&\int_{-\epsilon}^\epsilon \frac{ h u^2}{(h^2+u^2)^2}\, \dd u = -\frac{h\epsilon}{h^2 + \epsilon^2} + \arctan \left(\frac{\epsilon}{h}\right) \, , \\
&\int_{-\epsilon}^\epsilon \frac{ h^3}{(h^2+u^2)^2}\, \dd u = \frac{h\epsilon}{h^2 + \epsilon^2} + \arctan \left(\frac{\epsilon}{h}\right) \, , \\
&\int_{-\epsilon}^\epsilon \frac{hu^4}{(h^2+u^2)^3}\, \dd u = -\frac{3h^3\epsilon+5h\epsilon^3}{4(h^2 + \epsilon^2)^2} + \frac{3}{4} \arctan \left(\frac{\epsilon}{h}\right) \, , \\
&\int_{-\epsilon}^\epsilon \frac{h^3u^2}{(h^2+u^2)^3}\, \dd u = \frac{-h^3\epsilon+h\epsilon^3}{4(h^2 + \epsilon^2)^2} + \frac{1}{4} \arctan \left(\frac{\epsilon}{h}\right) \, , \\
&\int_{-\epsilon}^\epsilon \frac{h^5}{(h^2+u^2)^3}\, \dd u = \frac{5h^3\epsilon+3h\epsilon^3}{4(h^2 + \epsilon^2)^2} + \frac{3}{4} \arctan \left(\frac{\epsilon}{h}\right) \, , \\
&\int_{-\epsilon}^\epsilon \frac{h^7}{(h^2+u^2)^4}\, \dd u =  \frac{33 h^5 \epsilon +40 h^3 \epsilon ^3+15 h \epsilon ^5}{24 \left(h^2+\epsilon ^2\right)^3} + \frac{5}{8} \arctan\left(\frac{\epsilon}{h}\right) \, , \\
&\int_{-\epsilon}^\epsilon \frac{h^5 u^2}{(h^2+u^2)^4}\, \dd u =  \frac{-3 h^5 \epsilon +8 h^3 \epsilon ^3+3 h \epsilon ^5}{24 \left(h^2+\epsilon ^2\right)^3} + \frac{1}{8} \arctan\left(\frac{\epsilon}{h}\right) \, , \\
&\int_{-\epsilon}^\epsilon \frac{h^3 u^4}{(h^2+u^2)^4}\, \dd u =  \frac{-3 h^5 \epsilon -8 h^3 \epsilon ^3+3 h \epsilon ^5}{24 \left(h^2+\epsilon ^2\right)^3} + \frac{1}{8} \arctan\left(\frac{\epsilon}{h}\right) \, , \\
&\int_{-\epsilon}^\epsilon \frac{h u^6}{(h^2+u^2)^4}\, \dd u =  -\frac{  15 h^5 \epsilon +40 h^3 \epsilon ^3+33 h \epsilon ^5}{24 \left(h^2+\epsilon ^2\right)^3} + \frac{5}{8} \arctan\left(\frac{\epsilon }{h}\right) \, .
\end{align}
These integrals can be computed using a standard trigonometric substitution. It is easy to see that $\displaystyle\lim_{h \to 0^\pm} \arctan\left( \frac{\epsilon}{h} \right) = \pm \frac{\pi}{2}$ while the other terms in the expressions on the right go to zero as $h \to 0^\pm$ for fixed $\epsilon > 0$. Note that odd powers of $u$ integrate to zero while higher order terms can be bounded uniformly in $h$. Therefore, this list of integrals is exhaustive for the kernels in the integral equations, which have at most $\|\br\|^8$ in the denominator. 
\end{rmk}

\begin{rmk}It is straightforward to show the following bounds for $u \in [-\epsilon,\epsilon]$ uniformly as $h \to 0^\pm$: 
\begin{align}
    \left| \frac{u^{2j}}{(h^2 + u^2)^j} \right| & \leq 1  \, ,
\end{align}
for $j\in \mathbb{N}$. Hence it follows that the contribution of these terms is $\mathcal{O}(1)$ and vanishes (after integration) as $\epsilon$ goes to 0.  
\end{rmk}

The flexural wave Green's function $G(\bx,\by)$ only differs from the biharmonic Green's function $G^B(\bx,\by)$ by lower order terms, therefore it is sufficient to perform this analysis on $G^B(\bx,\by)$. Though these calculations are relatively standard, see e.g. \cite{kress, JIANG20117488}, we perform them here for completeness.

\subsection{Clamped plate}\label{app:jumpcondclamped}
Because $\log\| \bx-\by\|$ is integrable, the dominated convergence theorem
implies that logarithmically singular kernels have no singular contribution.
Thus, we can conclude the following limit for the clamped plate operator $\cK_2$ with kernel \eqref{clampedk2}:
\begin{align}
    \lim_{\bx \to \bx_0^{\pm}}\mathcal{K}_{2}[\sigma](\bx) &= \mathcal{K}_{12}[\sigma](\bx_0) \, .
\end{align}
Next, we turn our attention to the clamped plate kernels $K_{11}$ and $K_{22}$. Recall that we have the kernel:
\begin{align}
G^{B}_{\bn_\by  \bn_\by  \bn_\by } (\bx ,\by)  &= -\frac{3}{4\pi} \frac{[\br\cdot \bn (\by)]}{\|\br \|^2} + \frac{1}{2\pi} \frac{[\br\cdot \bn (\by)]^3}{\|\br \|^4} \, . \label{k11clamped1}
\end{align}
Using the expressions in Lemma \ref{lem:uexpansions}, we have to leading order:
\begin{align}
G^{B}_{\bn_\by  \bn_\by  \bn_\by } (\bgamma(0) + h \bn(0) ,\bgamma(u))  &= -\frac{3}{4\pi} \frac{h}{h^2+u^2} + \frac{1}{2\pi} \frac{h^3}{(h^2+u^2)^2} + \mathcal{O}(1) \, .
\end{align}
Using the integrals in Remark \ref{rmk:integrals} we find that:
\begin{equation}
\lim_{h \to 0^\pm}\int_{-\epsilon}^{\epsilon} G^{B}_{\bn_\by  \bn_\by  \bn_\by } (\bgamma(0) + h \bn(0) ,\bgamma(u)) \sigma(s(u)) s'(u) \, \dd u  = \mp\frac{1}{2} \sigma(0) + \mathcal{O}(\epsilon) \, . 
\label{jumpcond1}
\end{equation}
Note that this jump relation matches the limit of $\widehat{G_{\bn_\by \bn_\by \bn_\by}}$ as $\xi \to \pm \infty$ in the Fourier domain \eqref{fourier1}. The jump relation for $G^B_{\bn_\bx \bn_\by \bn_\by}$ in $K^B_{22}$ only differs by a sign and can be derived the same way. 
We turn our attention to the other term in $K^B_{11}$:
\begin{align}
    G^B_{\bn_\by \btau_\by \btau_\by}(\bx,\by) =  \frac{1}{2\pi} \frac{[\br\cdot  \bn (\by) ][\br\cdot \btau (\by)]^2}{\|\br \|^4} -\frac{1}{4\pi}\frac{\br\cdot  \bn (\by) }{\|\br \|^2} \, .
\end{align}
Again, using the expressions in Lemma \ref{lem:uexpansions} we have:
\begin{align}
    G^B_{\bn_\by \btau_\by \btau_\by}(\bgamma(0) + h \bn(0),\bgamma(u)) =  \frac{1}{2\pi} \frac{hu^2}{(h^2+u^2)^2} -\frac{1}{4\pi}\frac{ h }{h^2+u^2} \, + \mathcal{O}(1) \, .
\end{align}
Note that these two contributions cancel, giving:
\begin{align}
\lim_{h \to 0^\pm}\int_{-\epsilon}^{\epsilon} G^{B}_{\bn_\by  \btau_\by  \btau_\by } (\bgamma(0) + h \bn(0) ,\bgamma(u)) \sigma(s(u)) s'(u) \, \dd u  &=
\mathcal{O}(\epsilon) \, .
\end{align}
The jump relation for $G^B_{\bn_\bx \btau_\by \btau_\by}$  in $K^B_{22}$ is the same and can be obtained the same way. Again, the kernel of the integral on the right is proven to be continuous in Lemma \ref{lem:clampedsmooth}, so this integral will disappear as $\epsilon \to 0$. From this, we conclude that:
\begin{align}
    \lim_{\bx \to \bx_0^{\pm}}\mathcal{K}_{1}[\sigma](\bx) &= \mp \frac{1}{2} \sigma(\bx_0) + \mathcal{K}_{11}[\sigma](\bx_0) \, , \\
    \lim_{\bx \to \bx_0^{\pm}}\bn(\bx_0) \cdot \nabla_{\bx} \mathcal{K}_{2}[\sigma](\bx) &= \mp \frac{1}{2} \sigma(\bx_0) + \mathcal{K}_{22}[\sigma](\bx_0) \, .
\end{align}
Finally, we turn our attention to $ K^B_{21} = G^B_{\bn_\bx  \bn_\by  \bn_\by  \bn_\by } + 3G^B_{\bn_\bx  \bn_\by  \btau_\by  \btau_\by }$. Notice that this kernel can be written as:
\begin{align}
    K^B_{21}(\bx,\by) &=    \frac{3}{\pi} \frac{[\br\cdot \bn (\by)][\br\cdot \btau (\by)][\bn (\bx ) \cdot \btau (\by)]}{\|\br \|^4}  - \frac{3}{\pi} \frac{[\br\cdot \bn (\by)] [\br\cdot \bn (\bx )]}{\|\br \|^4} + \frac{4}{\pi} \frac{[\br\cdot \bn (\by)]^3 [\br\cdot \bn (\bx )]}{\|\br \|^6}  \, .
\end{align}
Again, using Lemma \ref{lem:uexpansions}, the expansion becomes:
\begin{align}
    K^B_{21}(\bgamma(0) + h \bn(0),\bgamma(u)) &= \frac{3}{\pi} \frac{h u^2  \kappa(0)}{(u^2+h^2)^2} - \frac{3}{\pi} \frac{h^2}{(h^2+u^2)^2} + \frac{4}{\pi} \frac{h^4 - h^3 u^2 \kappa(0)}{(h^2+u^2)^3}  + \mathcal{O}(1) \, .
\end{align}
Using the integrals in Remark \ref{rmk:integrals} and taking the limit yields:
\begin{align}
\lim_{h \to 0^\pm}\int_{-\epsilon}^{\epsilon} K^{B}_{21} (\bgamma(0) + h \bn(0) ,\bgamma(u)) \sigma(s(u)) s'(u) \, \dd u  &= \pm \kappa(0) \sigma(0) + \mathcal{O}(\epsilon) \, .
\end{align}
From this we write down the last jump relation:
\begin{align}
    \lim_{\bx \to \bx_0^{\pm}}\bn(\bx_0) \cdot \nabla_{\bx}  \mathcal{K}_{1}[\sigma](\bx) &=  \pm \kappa(\bx_0) \sigma(\bx_0) + \mathcal{K}_{21}[\sigma](\bx_0) \, , 
\end{align}
which concludes the derivation of jump relations for the interior and exterior clamped plate problems.

\subsection{Supported plate}\label{app:jumpcondsupp}

We first examine the supported plate kernel $K^B_{11} = G^B_{\bn_\by  \bn_\by  \bn_\by }   + \alpha_1 G^B_{\bn_\by  \btau_\by  \btau_\by }  +  \alpha_2 \kappa(\by ) G^B_{\bn_\by  \bn_\by  } + \alpha_3 \kappa'(\by ) G^B_{\btau_\by  } $. The terms $ G^B_{\bn_\by  \bn_\by  \bn_\by } $ and $ G^B_{\bn_\by  \btau_\by  \btau_\by }$ were analyzed in the previous section, and their jump relations are reported there. Meanwhile, the kernels $ \kappa(\by)  G^B_{\bn_\by  \bn_\by  } $ and $ \kappa'(\by ) G^B_{\btau_\by  }$ are log singular and continuous, respectively, and can be bounded uniformly for all $|h| < h_0$ provided that $\partial\Omega$ is sufficiently smooth. From this, we can immediately conclude the first jump relation:
\begin{align}
        \lim_{\bx \to \bx_0^{\pm}}\mathcal{K}_{1}[\sigma](\bx) &= \mp \frac{1}{2} \sigma(\bx_0) + \mathcal{K}_{11}[\sigma](\bx_0) \, .
\end{align}
Next, the kernel $K^B_{12} = G^B_{\bn_{\by }}$ is continuous so we can also conclude the jump relation:
\begin{align}
        \lim_{\bx \to \bx_0^{\pm}}\mathcal{K}_{2}[\sigma](\bx) &= \mathcal{K}_{12}[\sigma](\bx_0) \, .
\end{align}
Moreover, the jump relations for $ G^B_{\bn_\bx  \bn_\bx  \bn_{\by }} $ and $ G^B_{\btau_\bx  \btau_\bx  \bn_{\by }}$ in kernel $K^B_{22}$ are the same as the jump relations for $ G^B_{\bn_\by  \bn_\by  \bn_{\by }} $ and $ G^B_{\bn_{\by } \btau_\by  \btau_\by  }$ respectively and can be calculated the same way. We can immediately conclude that:
\begin{align}
        \lim_{\bx \to \bx_0^{\pm}}\left ( \nu \Delta + (1-\nu) \left ( \bn(\bx_0) \cdot \nabla_{\bx}\right)^2
    \right ) \mathcal{K}_{2}[\sigma](\bx) &= \mp \frac{1}{2} \sigma(\bx_0) + \mathcal{K}_{22}[\sigma](\bx_0) \, .
\end{align}

The only kernel we have left to analyze is $ K_{21}^B =  G^B_{\bn_\bx  \bn_\bx  \bn_\by  \bn_\by  \bn_\by } + \nu  G^B_{\btau_\bx  \btau_\bx  \bn_\by  \bn_\by  \bn_\by }   + \alpha_1 G^B_{ \bn_\bx  \bn_\bx  \bn_\by  \btau_\by  \btau_\by } + \nu  \alpha_1 G^B_{\btau_\bx  \btau_\bx  \bn_\by  \btau_\by  \btau_\by } +  \alpha_2 \kappa(\by ) G^B_{\bn_\bx  \bn_\bx  \bn_\by  \bn_\by  } +\nu \alpha_2 \kappa(\by ) G^B_{\btau_\bx  \btau_\bx  \bn_\by  \bn_\by  } + \alpha_3 \kappa'(\by ) G^B_{\bn_\bx  \bn_\bx  \btau_\by  } + \nu \alpha_3 \kappa'(\by ) G^B_{\btau_\bx  \btau_\bx  \btau_\by  }$. We write out the full expression for this kernel:

\begin{align}
    K_{21}^B(\bx,\by) &= \frac{\alpha_3 \kappa'(\by) (1-\nu) [\br\cdot\nx]^2 [\br \cdot \tauy]}{2 \pi  \|\br\|^4}-\frac{\alpha_3 \kappa'(\by) \nu [\br \cdot \taux] [\taux \cdot \tauy]}{2 \pi  \|\br\|^2} \nonumber \\
    &\qquad -\frac{\alpha_3 \kappa'(\by) (1-\nu) [\br \cdot \tauy]}{4 \pi  \|\br\|^2}-\frac{\alpha_3 \kappa'(\by) [\br \cdot \nx] [\nx \cdot \tauy]}{2 \pi  \|\br\|^2} \nonumber
    \\ 
    &\qquad + \frac{\alpha_2 \kappa(\by)}{2\pi} \frac{[\bn (\bx )\cdot \bn (\by)]^2}{\|\br \|^2} - \frac{2 \alpha_2 \kappa(\by)}{\pi} \frac{[\br\cdot \bn (\bx )][\br \cdot \bn (\by)][\bn (\bx )\cdot \bn (\by)]}{\|\br \|^4} \nonumber \\
    &\qquad -\frac{(1-\nu)\alpha_2 \kappa(\by)}{2\pi}\frac{[\br\cdot \bn (\by)]^2}{\|\br \|^4} + \frac{2(1-\nu)\alpha_2 \kappa(\by)}{\pi} \frac{[\br\cdot \bn (\bx )]^2[\br\cdot \bn (\by)]^2}{\|\br \|^6} \nonumber \\
    &\qquad- \frac{(1-\nu)\alpha_2 \kappa(\by)}{2\pi} \frac{[\br\cdot \bn (\bx )]^2}{\|\br \|^4} + \frac{(1-\nu)\alpha_2 \kappa(\by)}{4\pi} \frac{1}{\|\br \|^2}  \nonumber \\
    &\qquad + \frac{\nu \alpha_2 \kappa(\by)}{2\pi} \frac{[\btau (\bx)\cdot \bn (\by)]^2}{\|\br \|^2} - \frac{2\nu \alpha_2 \kappa(\by)}{\pi} \frac{[\br\cdot \btau (\bx)][\br \cdot \bn (\by)][\btau (\bx)\cdot \bn (\by)]}{\|\br \|^4}  \nonumber \\
&\qquad -\frac{12 (1-\alpha_1) \nu [\br \cdot \nx]^2 [\br \cdot \ny]^3}{\pi  \| \br \|^8}+\frac{(6-10 \alpha_1) \nu [\br \cdot \nx]^2 [\br \cdot \ny]}{\pi  \|\br\|^6} \nonumber \\
&\qquad -\frac{8 \alpha_1 \nu [\br \cdot \ny] [\br \cdot \taux] [\br \cdot \tauy] [\taux \cdot \tauy]}{\pi  \|\br\|^6} \nonumber \\
&\qquad -\frac{(12-4 \alpha_1) \nu [\br \cdot \ny]^2 [\taux \cdot \ny] [\br \cdot \taux]}{\pi  \|\br\|^6} \nonumber \\
&\qquad +\frac{10 (1-\alpha_1) \nu [\br \cdot \ny]^3}{\pi  \|\br\|^6}+\frac{\alpha_1 \nu [\br \cdot \ny] [\taux\cdot\tauy]^2}{\pi  \|\br\|^4} \nonumber \\
&\qquad +\frac{2 \alpha_1 \nu [\taux \cdot \ny] [\br \cdot \tauy] [\taux \cdot \tauy]}{\pi  \|\br\|^4}+\frac{3 (1-\alpha_1) \nu [\taux \cdot \ny] [\br \cdot \taux]}{\pi  \|\br\|^4} \nonumber \\
&\qquad -\frac{(9-17 \alpha_1) \nu [\br \cdot \ny]}{2 \pi  \|\br\|^4}+\frac{12 (1-\alpha_1) [\br \cdot \nx]^2 [\br \cdot \ny]^3}{\pi  \| \br \|^8} \nonumber \\
&\qquad -\frac{8 \alpha_1 [\br \cdot \nx] [\nx \cdot \tauy] [\br \cdot \ny] [\br \cdot \tauy]}{\pi  \|\br\|^6} \nonumber \\
&\qquad -\frac{(12-4 \alpha_1) [\nx \cdot \ny][\br \cdot \nx] [\br \cdot \ny]^2}{\pi  \|\br\|^6}-\frac{(6-10 \alpha_1) [\br \cdot \nx]^2 [\br \cdot \ny]}{\pi  \|\br\|^6} \nonumber \\
&\qquad +\frac{\alpha_1 ([\nx \cdot \tauy])^2 [\br \cdot \ny]}{\pi  \|\br\|^4}+\frac{2 \alpha_1 [\nx \cdot \ny][\nx \cdot \tauy] [\br \cdot \tauy]}{\pi  \|\br\|^4} \nonumber \\
&\qquad +\frac{3 (1-\alpha_1) [\nx \cdot \ny][\br \cdot \nx]}{\pi  \|\br\|^4}-\frac{2 (1-\alpha_1) [\br \cdot \ny]^3}{\pi  \|\br\|^6}+\frac{3 (1-\alpha_1) [\br \cdot \ny]}{2 \pi  \|\br\|^4} \nonumber \\
&\qquad +\frac{3 \nu [\br \cdot \ny] ([\taux \cdot \ny])^2}{\pi  \|\br\|^4}+\frac{3 [\nx \cdot \ny]^2 [\br \cdot \ny]}{\pi  \|\br\|^4} \, .
\end{align}

\begin{rmk}
    Recall that $\kappa(s)$ can be expanded as:
\begin{align}
    \kappa( s) &= \kappa(0) + s \kappa'(0) + \frac{1}{2} s^2 \kappa''(0) + \mathcal{O}(s^3) \, .
\end{align}
Performing the change of variables using \eqref{s(u)}, we have:
\begin{align}
    \kappa( s(u)) = \kappa(0) + u (1 - \frac12  h \kappa(0)  ) \kappa'(0) +  \mathcal{O}(u^2 + h^2) \, . \label{kappaexpansion}
\end{align}
Similarly we have:
\begin{align}
    \kappa'( s(u)) = \kappa'(0) + u (1 - \frac12  h \kappa(0)  ) \kappa''(0) +  \mathcal{O}(u^2 + h^2) \, . \label{kappapexpansion}
\end{align}
\end{rmk}

For succinctness, we multiply the kernel by $s'(u) \sigma(s(u))$ which appears in the integrand and expand directly with the help of Mathematica:
\begin{align}
    K_{21}^B &(\bgamma(0) + h \bn(0),\bgamma(u))  s'(u) \sigma(s(u)) \nonumber \\
    &= \frac{1}{8\pi} \frac{\alpha_2 (\nu-3) \left(h \kappa^2(0) \sigma(0)-2\kappa(0)(\sigma(0)+\sigma'(0) u)-2 \kappa'(0) \sigma(0) u\right)+2 \alpha_3 \kappa'(0) (\nu+1) \sigma(0) u}{h^2+u^2} \nonumber \\
    &\quad - \frac{1}{6 \pi } \frac{h^2 \kappa'(0) \sigma(0) u (\alpha_1 (19 \nu-9)-3 (\alpha_2 (4 \nu-6)+\alpha_3 \nu-\alpha_3+3 \nu-5))}{(h^2+u^2)^2} \nonumber \\
    &\quad + \frac{1}{16 \pi}\frac{h^3 \kappa^2(0) \sigma(0) (3 \alpha_1 (19 \nu-9)-8 \alpha_2 (2 \nu-3)-27 \nu+45)}{(h^2+u^2)^2 } \nonumber \\
    &\quad + \frac{1}{4 \pi}\frac{h^2\kappa(0)\sigma(0) (\alpha_1 (9-19 \nu)+4 \alpha_2 (2 \nu-3)+9 \nu-15)}{(h^2+u^2)^2 } \nonumber \\
    &\quad + \frac{1}{2 \pi}\frac{h^2\kappa(0)\sigma'(0) u (\alpha_1 (9-19 \nu)+\alpha_2 (4 \nu-6)+9 \nu-15)}{(h^2+u^2)^2 } \nonumber \\
    &\quad - \frac{1}{16 \pi}\frac{\alpha_1 h \sigma(0) \left(\nu \left(15 \kappa^2(0) u^2-152\right)-45 \kappa^2(0) u^2+72\right)}{(h^2+u^2)^2 } \nonumber \\
    &\quad +\frac{1}{4 \pi}\frac{h u (\alpha_1 (19 \nu-9)-9 \nu+15) (2 \sigma'(0)+\sigma''(0) u)}{(h^2+u^2)^2 } \nonumber \\
    &\quad -\frac{1}{12 \pi}\frac{u^2 (3\kappa(0)(15 \alpha_1 \nu-5 \alpha_1+3 \nu+3) (\sigma(0)+\sigma'(0) u)+4 \kappa'(0) \sigma(0) u (\alpha_1 (8 \nu-3)+3))}{(h^2+u^2)^2 } \nonumber \\
    &\quad +\frac{1}{16 \pi}\frac{h \sigma(0) \left(\nu \left((16 \alpha_2+93) \kappa^2(0) u^2-72\right)-75 \kappa^2(0) u^2+120\right)}{(h^2+u^2)^2 } \nonumber \\
&\quad + \frac{1}{2 \pi } \frac{h^5 \kappa^2(0) \sigma(0) (-3 \alpha_1 (5 \nu-4)+2 \alpha_2 (\nu-1)+3 (4 \nu-5))}{(h^2+u^2)^3 } \nonumber \\
&\quad + \frac{1}{\pi } \frac{2 h^4\kappa(0)\sigma(0) (\alpha_1 (5 \nu-4)+\alpha_2 (-\nu)+\alpha_2-4 \nu+5)}{(h^2+u^2)^3} \nonumber \\
&\quad + \frac{1}{\pi} \frac{2 h^4\kappa(0)\sigma'(0) u (2 \alpha_1 (5 \nu-4)+\alpha_2 (-\nu)+\alpha_2-8 \nu+10)}{(h^2+u^2)^3 } \nonumber \\
&\quad + \frac{1}{3 \pi} \frac{2 h^4 \kappa'(0) \sigma(0) u (2 \alpha_1 (5 \nu-4)-3 \alpha_2 (\nu-1)-8 \nu+10)}{(h^2+u^2)^3 } \nonumber \\
&\quad + \frac{1}{2 \pi} \frac{\alpha_1 h^3 \sigma(0) \left(\nu \left(11 \kappa^2(0) u^2-40\right)-16 \kappa^2(0) u^2+32\right)}{ (h^2+u^2)^3} \nonumber \\
&\quad - \frac{1}{\pi} \frac{2 h^3 u (\alpha_1 (5 \nu-4)-4 \nu+5) (2 \sigma'(0)+\sigma''(0) u)}{ (h^2+u^2)^3} \nonumber \\
&\quad + \frac{1}{\pi} \frac{2 h^2\kappa(0)u^2 (\sigma(0) (\alpha_1 (5 \nu-4)+3)+\sigma'(0) u (\alpha_1 (7 \nu-4)+3))}{(h^2+u^2)^3 } \nonumber \\
&\quad + \frac{1}{3 \pi} \frac{4 h^2 \kappa'(0) \sigma(0) u^3 (\alpha_1 (7 \nu-6)-3 \nu+6)}{(h^2+u^2)^3 } \nonumber 
\end{align}
\begin{align}
&\quad - \frac{1}{\pi}\frac{2 \alpha_1 h u^3 \left(\kappa^2(0) \sigma(0) u+4 \nu \sigma'(0)+2 \nu \sigma''(0) u\right)}{(h^2+u^2)^3 }  + \frac{1}{\pi} \frac{2 \alpha_1 h \nu \sigma(0) u^2 \left(\kappa^2(0) u^2-4\right)}{(h^2+u^2)^3 } \nonumber \\
&\quad + \frac{1}{3 \pi} \frac{4 \alpha_1 \nu u^4 (3\kappa(0)(\sigma(0)+\sigma'(0) u)+2 \kappa'(0) \sigma(0) u)}{(h^2+u^2)^3 } \nonumber \\
&\quad + \frac{1}{2 \pi } \frac{h^3 \sigma(0) \left(-32 \kappa^2(0) \nu u^2+25 \kappa^2(0) u^2+32 \nu-40\right)}{(h^2+u^2)^3}  + \frac{1}{\pi} \frac{3 h \kappa^2(0) (1-2 \nu) \sigma(0) u^4}{(h^2+u^2)^3 } \nonumber
\\    &\quad -\frac{1}{\pi}\frac{2 (\alpha_1-1) h^6 (\nu-1) (3\kappa(0)(\sigma(0)+2 \sigma'(0) u)+2 \kappa'(0) \sigma(0) u)}{(h^2+u^2)^4 } \nonumber \\
    &\quad +\frac{1}{2 \pi}\frac{3 (\alpha_1-1) h^5 (\nu-1) \left(\sigma(0) \left(8-5 \kappa^2(0) u^2\right)+4 u (2 \sigma'(0)+\sigma''(0) u)\right)}{(h^2+u^2)^4 } \nonumber \\
    &\quad - \frac{1}{\pi}\frac{2 (\alpha_1-1) h^4 (\nu-1) u^2 (3\kappa(0)(\sigma(0)+\sigma'(0) u)+4 \kappa'(0) \sigma(0) u)}{ (h^2+u^2)^4} \nonumber \\
    &\quad +\frac{1}{2 \pi }\frac{3 (\alpha_1-1) h^3 \kappa^2(0) (\nu-1) \sigma(0) \left(3 h^4-4 u^4\right)}{(h^2+u^2)^4} + \mathcal{O}(1) \, . 
\end{align}
As before, these expressions are readily integrated using Remark \ref{rmk:integrals}, yielding:
\begin{align}
    \lim_{h \to 0^\pm}\int_{-\epsilon}^{\epsilon} K_{21}^B &(\bgamma(0) + h \bn(0),\bgamma(u))  s'(u) \sigma(s(u))= \pm c_0 \kappa^2(0) \sigma(0)+ 
    \mathcal{O}(\epsilon) \, ,
\end{align}
where $c_0 = \frac{1}{2}(-1 + \alpha_1 - \nu + \alpha_1 \nu + \alpha_2 \nu)$. From this we obtain the final jump relation:
\begin{align}
    \lim_{\bx \to \bx_0^{\pm}}\left ( \nu \Delta + (1-\nu) \left (\bn(\bx_0) \cdot \nabla_{\bx} \right )^2
    \right ) \mathcal{K}_{1}[\sigma](\bx) &=  \pm c_0 \kappa^2(\bx_0) \sigma(\bx_0) + \mathcal{K}_{21}[\sigma](\bx_0) \, ,
\end{align}
concluding the derivation of jump relations for the supported plate kernels.

\subsection{Free plate}\label{app:jumpcondfree}

The $K_{12}$ kernel is logarithmically singular, so that
\begin{align}
    \lim_{\bx \to \bx_0^{\pm}} \mathcal{B}_1(\bx_0) \mathcal{K}_{2}[\sigma](\bx) &= \mathcal{K}_{12}[\sigma](\bx_0) \, , 
\end{align}
where $\bx_0^+$ represents the exterior limit and $\bx_0^-$ represents the interior limit, each along the normal direction. Moreover, we have already examined kernels related to $G^B_{\bn_\bx \bn_\bx \bn_\bx}$ and $ G^B_{\bn_\bx \btau_\bx \btau_\bx}$. Therefore, we can also state the following jump relation:
\begin{align}
    \lim_{\bx \to \bx_0^{\pm}} \mathcal{B}_2(\bx_0) \mathcal{K}_{2}[\sigma](\bx) &= \pm \frac{1}{2} \sigma(\bx_0) + \mathcal{K}_{22}[\sigma](\bx_0) \, .
\end{align}

We observe that $K_{11}^a$ is the same as the $K_{22}$ kernel for the supported plate and has the same 
jump properties.
It is then necessary for us to analyze $K_{11}^b = G_{\bn_\bx \bn_\bx \btau_\by} + \nu G_{\btau_\bx \btau_\bx \btau_\by}$. First we have:
\begin{align}
        G^{B}_{\bn_\bx  \bn_\bx  \btau_\by } (\bx,\by) &= -\frac{1}{2\pi} \frac{[\br\cdot  \bn (\bx ) ][ \bn (\bx )  \cdot  \btau (\by) ]}{\|\br \|^2} + \frac{1}{2\pi} \frac{[\br\cdot  \btau (\by) ][\br\cdot  \bn (\bx ) ]^2}{\|\br \|^4} -\frac{1}{4\pi}\frac{\br\cdot  \btau (\by) }{\|\br \|^2} \, .
\end{align}
Using Lemma \ref{lem:uexpansions}, this becomes:
\begin{align}
        G^{B}_{\bn_\bx  \bn_\bx  \btau_\by } (\bgamma(0) + h \bn(0),\bgamma(u)) &= - \frac{1}{2\pi} \frac{u h^2 }{(h^2 +  u^2)^2} +\frac{1}{4\pi} \frac{u}{h^2+u^2} + \mathcal{O}(1) \, .
\end{align}
Note that these terms are odd so they integrate to zero. Moreover, when they are multiplied by the next order term in $\sigma$, the resulting terms are uniformly bounded. We conclude the following:
\begin{align}
\lim_{h \to 0^\pm}\int_{-\epsilon}^{\epsilon} G^{B}_{\bn_\bx  \bn_\bx  \btau_\by } (\bgamma(0) + h \bn(0) ,\bgamma(u)) \sigma(s(u)) s'(u) \, \dd u  &= \mathcal{O}(\epsilon) \, .
\end{align}
Similarly, we have:
\begin{align}
         G^{B}_{\btau_\bx  \btau_\bx  \btau_\by } (\bx,\by) &= -\frac{1}{2\pi} \frac{[\br\cdot \btau (\bx)][\btau (\bx) \cdot  \btau (\by) ]}{\|\br \|^2} + \frac{1}{2\pi} \frac{[\br\cdot  \btau (\by) ][\br\cdot \btau (\bx)]^2}{\|\br \|^4} -\frac{1}{4\pi}\frac{\br\cdot  \btau (\by) }{\|\br \|^2} \, .
\end{align}
Expanding, we get:
\begin{align}
         G^{B}_{\btau_\bx  \btau_\bx  \btau_\by } (\bgamma(0) + h \bn(0),\bgamma(u)) &= \frac{3}{4\pi} \frac{u}{h^2 + u^2} - \frac{1}{2\pi} \frac{u^3}{(h^2+u^2)^2} + \mathcal{O}(1) \, .
\end{align}
Again, these terms are odd and the next order terms are uniformly bounded. Therefore:
\begin{align}
\lim_{h \to 0^\pm}\int_{-\epsilon}^{\epsilon} G^{B}_{\btau_\bx  \btau_\bx  \btau_\by } (\bgamma(0) + h \bn(0) ,\bgamma(u)) \sigma(s(u)) s'(u) \, \dd u  &= 
\mathcal{O}(\epsilon) \, .
\end{align}
We conclude the following jump relation for $\mathcal{B}_1$ applied to $\mathcal{K}_{1}$:
\begin{align}
    \lim_{\bx \to \bx_0^{\pm}} \mathcal{B}_1(\bx_0) \mathcal{K}_{1}[\sigma](\bx) &= \mp \frac{1}{2} \sigma(\bx_0) + \mathcal{K}_{11}[\sigma](\bx_0) \, .
\end{align}

It remains to analyze kernels $K_{21}^a$ and $K_{21}^b$. First we focus on $K_{21}^b$. We have already analyzed $G^B_{\bn_\bx \bn_\bx \btau_\by}$ and $G^B_{\btau_\bx \btau_\bx \btau_\by}$, so we only need to look at $G^B_{\bn_\bx  \bn_\bx  \bn_\bx \btau_\by}$ and $ G^B_{\bn_\bx  \btau_\bx  \btau_\bx \btau_\by }$. Recall the formulae for these kernels are:
\begin{align}
     G^{B}_{\bn_\bx  \bn_\bx  \bn_\bx  \btau_\by } {(\bx,\by)} &= -\frac{3}{4\pi}\frac{ \bn (\bx )  \cdot  \btau (\by) }{\|\br \|^2} + \frac{3}{2\pi} \frac{[\br\cdot  \bn (\bx ) ][\br\cdot  \btau (\by) ]}{\|\br \|^4} + \frac{3}{2\pi} \frac{[\br\cdot  \bn (\bx ) ]^2[  \bn (\bx )  \cdot  \btau (\by) ]}{\|\br \|^4} \nonumber \\
    &\qquad -\frac{2}{\pi}\frac{[\br\cdot  \btau (\by)  ][\br \cdot  \bn (\bx ) ]^3}{\|\br \|^6} \, ,\\
     G^{B}_{\bn_\bx  \btau_\bx  \btau_\bx  \btau_\by  } {(\bx,\by)} &= \frac{1}{\pi} \frac{[\br\cdot  \bn (\bx )  ][\br\cdot \btau (\bx)] [\btau (\bx)\cdot  \btau (\by) ]}{\|\br \|^4}+ \frac{1}{2\pi }\frac{[ \btau (\by)  \cdot  \bn (\bx ) ] [\br \cdot \btau (\bx)]^2}{\|\br \|^4} - \frac{1}{4\pi} \frac{ \btau (\by)  \cdot  \bn (\bx ) }{\|\br \|^2}  \nonumber \\
    &\qquad - \frac{2}{\pi} \frac{[\br\cdot  \bn (\bx ) ][\br \cdot  \btau (\by) ][\br \cdot \btau (\bx)]^2}{\|\br \|^6} + \frac{1}{2\pi} \frac{[\br \cdot  \bn (\bx ) ] [\br \cdot  \btau (\by)  ]}{\|\br \|^4}  \, .
\end{align}
Expanding the product of these kernels with $\sigma(s(u)) s'(u)$ gives:
\begin{align}
     G^{B}_{\bn_\bx  \bn_\bx  \bn_\bx  \btau_\by } (\bgamma(0) + h \bn(0),\bgamma(u))& \sigma(s(u)) s'(u) \nonumber \\
     &= -\frac{3}{4\pi}\frac{ -\kappa(0) \sigma(0) u }{h^2+u^2} - \frac{3}{2\pi} \frac{h \sigma(0) u + h \sigma'(0) u^2 + \frac12 \kappa(0) \sigma(0) u^3  }{(h^2+u^2)^2} \nonumber \\
     &\qquad - \frac{3}{2\pi} \frac{h^2 \kappa(0) \sigma(0) u}{(h^2+u^2)^2}   +\frac{2}{\pi}\frac{h^3 \sigma u + h^3 \sigma'(0) u^2 + \frac32 h^2 \kappa(0) \sigma u^3}{(h^2+u^2)^3}  + \mathcal{O}(1) \, ,\\
     G^{B}_{\bn_\bx  \btau_\bx  \btau_\bx  \btau_\by  } (\bgamma(0) + h \bn(0),\bgamma(u)) & \sigma(s(u)) s'(u) \nonumber \\
     &= \frac{1}{\pi} \frac{-h \sigma(0) u + h^2 \kappa(0) \sigma(0) u - h \sigma'(0) u^2 - 1/2 \kappa(0) \sigma(0) u^3 }{(h^2+u^2)^2} \nonumber \\
     &\qquad - \frac{1}{2\pi }\frac{\kappa(0) \sigma(0) u^3}{(h^2+u^2)^2} + \frac{1}{4\pi} \frac{ \kappa(0) \sigma(0) u }{h^2+u^2}  \nonumber \\
    &\qquad - \frac{2}{\pi} \frac{-h \sigma(0) u^3 + h^2 \kappa(0) \sigma(0) u^3 - h \sigma'(0) u^4 - 1/2 \kappa(0) \sigma(0) u^5}{(h^2+u^2)^3} \nonumber \\ 
    &\qquad + \frac{1}{2\pi} \frac{-h \sigma(0) u - h \sigma'(0) u^2 - 1/2 \kappa(0) \sigma(0) u^3}{(h^2+u^2)^2} + \mathcal{O}(1) \, .
\end{align}
Integrating these terms gives the following jump relations:
\begin{align}
    \lim_{h \to 0^\pm}\int_{-\epsilon}^{\epsilon} G^{B}_{\bn_\bx  \bn_\bx  \bn_\bx  \btau_\by } (\bgamma(0) + h \bn(0),\bgamma(u)) \sigma(s(u)) s'(u) \, \dd u 
    &= \mp \frac{1}{2}\sigma'(0) + \mathcal{O}(\epsilon) \, , \\ 
    \lim_{h \to 0^\pm}\int_{-\epsilon}^{\epsilon} G^{B}_{\bn_\bx  \btau_\bx  \btau_\bx  \btau_\by  } (\bgamma(0) + h \bn(0),\bgamma(u)) \sigma(s(u)) s'(u) \, \dd u 
    &= \mathcal{O}(\epsilon) \, .
\end{align}
Replacing $\sigma(s(u))$ with its Hilbert transform, we get:
\begin{align}
        \lim_{h \to 0^\pm}\int_{-\epsilon}^{\epsilon} K^b_{21} (\bgamma(0) + h \bn(0),\bgamma(u)) \mathcal{H}[\sigma](u) s'(u) \, \dd u 
        &= \mp \frac{1}{2}\frac{\dd}{\dd s} \mathcal{H}[\sigma](s(u)) \big|_{s=0} + \mathcal{O}(\epsilon) \, . 
        \label{k21bjump}
\end{align}

Next, we look at $K^a_{21} = G_{\bn_\bx  \bn_\bx  \bn_\bx \bn_\by} + (2-\nu) G_{\bn_\bx  \btau_\bx  \btau_\bx  \bn_\by} + (1-\nu)\kappa(\bx) \left(G_{\btau_\bx  \btau_\bx \bn_\by} -  G_{\bn_\bx  \bn_\bx \bn_\by}\right)$. From the previous section we know that $(1-\nu) \kappa(\bx) G_{\bn_\bx  \bn_\bx \bn_\by}$ contributes a factor of $\mp \displaystyle\frac{1-\nu}{2} \kappa(\bx) \sigma(\bx)$ to the jump relation. Next, we analyze $G^B_{\bn_\bx  \bn_\bx  \bn_\bx \bn_\by} $ and $ G^B_{\bn_\bx  \btau_\bx  \btau_\bx  \bn_\by}$. These kernels are given by:
\begin{align}
    G^{B}_{\bn_\bx  \bn_\bx  \bn_\bx  \bn_\by } (\bx,\by) &= -\frac{3}{4\pi} \frac{ \bn (\bx ) \cdot  \bn (\by) }{\|\br \|^2} + \frac{3}{2\pi} \frac{[\br\cdot  \bn (\by) ][\br\cdot  \bn (\bx ) ]}{\|\br \|^4}+ \frac{3}{2\pi} \frac{[\br\cdot  \bn (\bx ) ]^2 [ \bn (\bx )  \cdot  \bn (\by) ]}{\|\br \|^4} \nonumber \\
    &\qquad -\frac{2}{\pi} \frac{[\br \cdot \bn (\bx )]^3[\br\cdot \bn (\by)]}{\|\br \|^6} \, , \\
    G^{B}_{\bn_\bx  \btau_\bx  \btau_\bx   \bn_\by } (\bx,\by) &= \frac{1}{2\pi} \frac{[\br \cdot \btau (\bx)]^2 [ \bn (\bx )  \cdot  \bn (\by) ]}{\|\br \|^4} + \frac{1}{\pi} \frac{[\br\cdot \btau (\bx) ][\btau (\bx)\cdot  \bn (\by) ] [\br \cdot  \bn (\bx ) ] }{\|\br \|^4} \nonumber \\
    &\qquad - \frac{2}{\pi} \frac{[\br \cdot  \bn (\by) ][\br \cdot  \bn (\bx ) ] [\br \cdot \btau (\bx)]^2 }{\|\br \|^6} -\frac{1}{4\pi} \frac{ \bn (\bx )  \cdot  \bn (\by) }{\|\br \|^2 } \nonumber \\ &\qquad+ \frac{1}{2\pi} \frac{[\br \cdot  \bn (\by) ][\br \cdot  \bn (\bx ) ]}{\|\br \|^4} \, .
\end{align}
These kernels are hyper-singular, and must be treated using a different
framework from the principal value approach introduced in this section.
For the sake of brevity,
we appeal to the following well-known fact~\cite{hsiao2008boundary,hackbusch2012integral}:
\begin{lem} \label{dprimejump}
  Let $\mathcal{D}$ be a double layer potential for the Laplace
  equation defined on a sufficiently smooth boundary curve $\partial \Omega$.
  For $\bx_0\in \partial \Omega$, we have
   \begin{align}
   \lim_{\bx\to\bx_0^\pm} \bn(\bx_0) \cdot \nabla \mathcal{D}(\bx)
  = \fp \int_{\partial \Omega} \left( \frac{\bn(\bx_0)\cdot\bn(\by)}{2\pi \|\bx_0-\by\|^2}
  - \frac{[\bn(\bx_0)\cdot(\bx_0-\by)][\bn(\by)\cdot(\bx_0-\by)]}{\pi \|\bx_0-\by\|^4}
  \right) \, \dd s(\by) \, ,
  \end{align}
  where the notation $\fp$ indicates that the integral is to be understood
  in the Hadamard finite part sense~\cite{hsiao2008boundary,hackbusch2012integral}.
\end{lem}
\begin{rmk}
  As observed in Section~\ref{freesection}, the fact that
  the integrals in $\mathcal{K}_{21}^a$ must be interpreted in
  the finite part sense does not affect the computation. This
  is because there is a cancellation between $\mathcal{K}_{21}^a$
  and the $\mathcal{H}'$ operator obtained from the jump in
  $K_{21}^b$. In fact, when defining $\mathcal{H}'$ as convolution against
  $K^{\mathcal H'}$, the integral must also be understood as a finite 
  part integral. When these kernels are combined, the result is
  continuous, cf. Lemma~\ref{freesmooth}, and can be treated
  using standard methods. 
\end{rmk}

Let $N(\bx,\by)$ be the hyper-singular integral kernel appearing above,
i.e.
\begin{align}
  N(\bx,\by) = \frac{\bn(\bx)\cdot\bn(\by)}{2\pi \|\bx-\by\|^2}
- \frac{[\bn(\bx)\cdot(\bx-\by)][\bn(\by)\cdot(\bx-\by)]}{\pi \|\bx-\by\|^4} \, .  
\end{align}
We then have
\begin{align}
  G^{B}_{\bn_\bx  \bn_\bx  \bn_\bx  \bn_\by } (\bx,\by) &= -\frac{3}{2} N(\bx,\by)
  + \frac{3}{2\pi} \frac{[\br\cdot  \bn (\bx ) ]^2 [ \bn (\bx )  \cdot  \bn (\by) ]}{\|\br \|^4} 
 -\frac{2}{\pi} \frac{[\br \cdot \bn (\bx )]^3[\br\cdot \bn (\by)]}{\|\br \|^6} \, , \\
  G^{B}_{\bn_\bx  \btau_\bx  \btau_\bx   \bn_\by } (\bx,\by) &=
  \frac{1}{2} N(\bx,\by) + \frac{1}{\pi} \frac{[\br\cdot \btau (\bx) ][\btau (\bx)\cdot  \bn (\by) ] [\br \cdot  \bn (\bx ) ] }{\|\br \|^4} 
 - \frac{1}{2\pi} \frac{[\br\cdot  \bn (\bx ) ]^2 [ \bn (\bx )  \cdot  \bn (\by) ]}{\|\br \|^4} \nonumber \\
  &\qquad + \frac{1}{\pi} \frac{[\br \cdot  \bn (\by) ][\br \cdot  \bn (\bx ) ] ([\br \cdot \bn (\bx)]^2 - [\br \cdot \btau (\bx)]^2) }{\|\br \|^6}  \, .
\end{align}
Because the hyper-singular parts will not have a jump by the lemma
above, we can focus on the parts which have principal value integrals, i.e.
\begin{align}
  G^{B,PV}_{\bn_\bx  \bn_\bx  \bn_\bx  \bn_\by } (\bx,\by) &=
     \frac{3}{2\pi} \frac{[\br\cdot  \bn (\bx ) ]^2 [ \bn (\bx )  \cdot  \bn (\by) ]}{\|\br \|^4} -\frac{2}{\pi} \frac{[\br \cdot \bn (\bx )]^3[\br\cdot \bn (\by)]}{\|\br \|^6} \, , \\
    G^{B,PV}_{\bn_\bx  \btau_\bx  \btau_\bx   \bn_\by } (\bx,\by) &=
    \frac{1}{\pi} \frac{[\br\cdot \btau (\bx) ][\btau (\bx)\cdot  \bn (\by) ] [\br \cdot  \bn (\bx ) ] }{\|\br \|^4} 
 - \frac{1}{2\pi} \frac{[\br\cdot  \bn (\bx ) ]^2 [ \bn (\bx )  \cdot  \bn (\by) ]}{\|\br \|^4} \nonumber \\
    &\qquad + \frac{1}{\pi} \frac{[\br \cdot  \bn (\by) ][\br \cdot  \bn (\bx ) ] ([\br \cdot \bn (\bx)]^2 - [\br \cdot \btau (\bx)]^2) }{\|\br \|^6}   \, .
\end{align}
Expanding these kernels gives:
\begin{align}
    G^{B,PV}_{\bn_\bx  \bn_\bx  \bn_\bx  \bn_\by }  (\bgamma(0)+h \bn(0), \bgamma(u)) 
    &= \frac{3}{2\pi} \frac{h^2 + h u^2 \kappa(0) }{(h^2+u^2)^2} - \frac{2}{\pi} \frac{h^4 + h^3 u^2 \kappa(0) }{(h^2+u^2)^3} + \mathcal{O}(1) \, , \\
    G^{B,PV}_{\bn_\bx  \btau_\bx  \btau_\bx   \bn_\by } (\bgamma(0) + h \bn(0), \bgamma(u))  &= \frac{1}{\pi} \frac{-hu^2 \kappa(0) }{(h^2+u^2)^2} - \frac{1}{2\pi} \frac{h^2 + hu^2 \kappa(0)}{(h^2+u^2)^2} \nonumber \\
    &\quad + \frac{1}{\pi} \frac{h^4 - h^2 u^2 + 2h^3 u^2 \kappa(0)}{(h^2+u^2)^3} + \mathcal{O}(1) \, .
\end{align}
These imply the following jump relations:
\begin{align}
     \lim_{h \to 0^\pm} \int_{-\epsilon}^{\epsilon}  G^{B,PV}_{\bn_\bx  \bn_\bx  \bn_\bx  \bn_\by } (\bgamma(0)+h \bn(0), \bgamma(u))  \sigma(s(u)) s'(u) \, \dd u 
    &= \pm \frac{1}{2}  \kappa(0) \sigma(0) + \mathcal{O}(\epsilon) \, , \\
     \lim_{h \to 0^\pm} \int_{-\epsilon}^{\epsilon}  G^{B,PV}_{\bn_\bx  \btau_\bx  \btau_\bx   \bn_\by } (\bgamma(0)+h \bn(0), \bgamma(u))  \sigma(s(u)) s'(u) \, \dd u 
    &= \mp \frac{1}{2}  \kappa(0) \sigma(0) + \mathcal{O}(\epsilon) \, .
\end{align}
Combining these two jump relations with the result from Lemma \ref{dprimejump} implies the jump relation for the combined term $K^{PV}:=G^B_{\bn_\bx  \bn_\bx  \bn_\bx \bn_\by} + (2-\nu) G^B_{\bn_\bx  \btau_\bx  \btau_\bx  \bn_\by}$ is given by:
\begin{align}
     \lim_{h \to 0^\pm} \int_{-\epsilon}^{\epsilon}  K^{PV} (\bgamma(0)+h \bn(0), \bgamma(u))  \sigma(s(u)) s'(u) \, \dd u 
    &= \pm \frac{\nu-1}{2}  \kappa(0) \sigma(0) + \mathcal{O}(\epsilon) \, ,
\end{align}
Note that this jump cancels with the jump from $(1-\nu) \kappa(\bx) G^B_{\bn_\bx  \bn_\bx \bn_\by}$, leading to the final jump relation:
\begin{align}
        \lim_{\bx \to \bx_0^{\pm}} \mathcal{B}_2(\bx_0) \mathcal{K}_{1}[\sigma](\bx) &= -\frac{\beta}{2} \frac{\dd}{\dd s} \mathcal{H}[\sigma](\bx_0) + \mathcal{K}_{21}[\sigma](\bx_0) \, ,
\end{align}
which concludes the derivation of jump relations for the free plate kernels.
}

\pagebreak 
\bibliography{elsrefs}   
\bibliographystyle{elsarticle-num}  

\end{document}